\documentclass[A4paper,bezier]{article}

\usepackage[tbtags]{amsmath}
\usepackage{amsfonts}
\usepackage{amssymb}
\usepackage{amsthm}
\usepackage{amscd}
\usepackage{mathrsfs}
\usepackage{newlfont}
\usepackage{txfonts} 
\usepackage{graphicx}

\title{Discrminantal Groups and Zariski Pairs of Sextic Curves}

\author{Jin-Gen Yang \thanks{Partially supported by NSF of China. \newline 2000 Mathematical
Subject Classification: Primary 14F45, Secondary 14H50,14J28.
}, Fudan University, Shanghai, jgyang@fudan.edu.cn\\
Jinjing Xie, Fudan University, Shanghai, 032018008@fudan.edu.cn}

\date{}

\newtheorem{lemma}{Lemma}[section]

\newtheorem{thm}[lemma]{Theorem}

\begin{document}

\hsize =14.5truecm

\large

 \pagenumbering{arabic}

\maketitle

\begin{abstract}
A series of Zariski pairs and four Zariski triplets were found by
using lattice theory of K3 surfaces. There is a Zariski triplet of
which one member is a deformation of another.
\end{abstract}

\section{Introduction}

In \cite{zariski1} Zariski showed that there are two irreducible
sextic curves $C_1,C_2$ with six cusps and the fundamental groups of
${\mathbb P}^2 \backslash C_1$ and ${\mathbb P}^2 \backslash C_2$
are not isomorphic. Such pairs are called Zariski pairs. The precise
definition of Zariski pair differs from paper to paper. Here we
adopt the following definition: Two plane curves $C_1,C_2$ of the
same degree form a Zariski pairs if $C_1,C_2$ have the same
combinatorial data (cf. \cite{Bartolo2}) and $({\mathbb P}^2,C_1)$
and $({\mathbb P}^2, C_2)$ are not homeomorphic. The Zariski triplet
and $k$-plet are defined similarly (cf. \cite{kplets}). A brief
account of the history of Zariski pairs can be found in
\cite{kplets}. It is remarkable that the degrees of all known
Zariski pairs are at least six.

Let $C$ be a reduced sextic curve with simple singularities only and
let $X$ be the K3 surface obtained from the double cover branched
over $C.$ Let $N_C$ be the orthogonal complement in $H^2(X,{\mathbb
Z})$ of the sublattice generated by all irreducible components of
the inverse image of $C$ in $X.$ Shimada shows in \cite{Shimada2}
that $N_C$ is a topological invariant of the pair $({\mathbb
P}^2,C).$ When $C$ is maximizing , i.e., the Milnor number of $C$ is
$19,$ $N_C$ is the transcendental lattice of the K3 surface $X.$ Let
$\gamma_X$ be the discriminantal form of the Picard lattice of $X.$
For some special maximizing sextics there are two non-isomorphic
positive definite lattices of rank two whose discriminantal forms
are isomorphic to $-\gamma_X.$ By Shimada's theorem they are Zariski
pairs, called arithmetic Zariski pairs. Shimada was able to
enumerate all such pairs (\cite{Shimada1,Shimada2}).

For any reduced sextic with simple singularities, not necessarily
maximizing, let $M$ be the primitive hull of the sublattice
generated by all irreducible components of the inverse image of $C$
in $X.$ By Shimada's theorem and Nikulin's lattice theory, the
discriminantal group $A$ of $M$ is a topological invariant of
$({\mathbb P}^2,C),$ which is weaker than $N_C.$ In this paper we
use this invariant to obtain a series of Zariski pairs and four
Zariski triplets of reduced sextics. Among them the most interesting
one is a Zariski triplet of three conics with $3A_5+3A_1,$ of which
one member of the triplet is the deformation of another member
(Theorem \ref{special}). To our knowledge this is the first such
example.

One significant difference between our Zariski pairs and Shimada's
arithmetic Zariski pairs is that in our examples although two
members of a pair have the same combinatorial data but for one
member there is a plane curve of low degree whose intersection
number with the sextic at every point is even. This geometric
property is not shared by arithmetic Zariski pairs of Milnor number
$19$, since the Picard groups of both members of such a pair are
isomorphic.

We would like to thank Professor W. Barth, Professor K. Zuo and Bo
Wu for helpful discussions on the subject.

After the finishing of this paper, the authors were kindly informed
by Professor Shimada that he obtained similar results, (cf.
\cite{Shimada09}).

\section{Discriminantal group of a sextic curve with simple
singularities}

Let $C$ be a reduced sextic curve with simple singularities only.
Let $p: Y \rightarrow {\mathbb P}^2$ be the double cover branched
over $C$ and let $\mu: X \rightarrow Y$ be the minimal resolution of
singularities of $Y.$ Then $X$ is a K3 surface and $H^2(X,{\mathbb
Z})$ is a unimodular lattice of signature $(3,19).$ Let
$\text{Pic}(X)$ denote the Picard lattice of $X.$ It is a primitive
sublattice of $H^2(X,{\mathbb Z}).$ Let $G$ be the sublattice of
$\text{Pic}(X)$ generated by all irreducible components of the
pull-back of $C$ in $X$ and let $\tilde{G}$ be the primitive hull of
$G$ in $H^2(X,{\mathbb Z}).$ We define the {\bf discriminantal
group} of $C$ to be the finite group $\tilde{G}^\vee/\tilde{G},$
where $\tilde{G}^\vee$ is the dual lattice of $\tilde{G}.$

\begin{lemma}\label{invariant}
Let $C$ and $C_1$ be two reduced sextic curves with simple
singularities only. If $({\mathbb P}^2,C)$ is homeomorphic to
$({\mathbb P}^2,C_1),$ then the discriminantal groups of $C$ and
$C_1$ are isomorphic.
\end{lemma}

Proof. Let $G$ be the sublattice generated by all irreducible
components of $(p \mu)^{-1}(C)$ and let $\tilde{G}$ be the primitive
hull of $G$ in $H^2(X,{\mathbb Z}).$ According to a theorem of
Shimada (\cite{Shimada2}) the orthogonal complement $G^\perp$ of $G$
in $H^2(X,{\mathbb Z})$ is a topological invariant of the pair
$({\mathbb P}^2,C).$ By Nikulin's lattice theory (\cite{nik}) the
 discriminantal group of $G^\perp$ is isomorphic to
that of $\tilde{G}.$ Hence the lemma holds. $\Box$

\medskip

Denote the unimodular even lattice of signature $(3,19)$ by
$\Lambda,$ called the K3 lattice. Recall that for any lattice $L,$
an overlattice $M$ of $L$ is a sublattice of the dual lattice
$L^\vee$ such that $M \supset L$ and $|M/L|< \infty.$

\begin{thm}[Urabe \cite{sextic,ams}]\label{urabe}
Let $G=\sum_k a_k A_k + \sum_l d_l D_l + \sum_m e_m E_m$ be a finite
Dynkin graph. Let $L(G)$ denote the negative definite lattice of
$G.$ Let ${\mathbb Z}\lambda$ be a lattice of rank one generated by
$\lambda$ with $\lambda^2=2.$ Then there is a reduced sextic curve
in ${\mathbb P}^2$ whose singularities correspond to $G$ if and only
if there is an overlattice $M$ of ${\mathbb Z}\lambda \oplus L(G)$
such that there is a primitive embedding of $M$ into the K3 lattice
$\Lambda$ such that

i) if $u \in M, u \lambda =0, u^2=-2,$ then $u \in L(G);$

ii) there is no $u \in M$ with $u \lambda =1$ and $u^2=0.$
\end{thm}

The main tool in Urabe's proof of the theorem is the surjectivity of
the period map for K3 surfaces. The sextic curve in the theorem can
be so chosen that the overlattice $M$ is exactly the Picard group of
the corresponding K3 surface and the pull-back of a line on the
plane belongs to the divisor class $\lambda.$

\section{Classical Zariski pair}

In this section we give a lattice-theoretic interpretation of
Zariski's classical example mentioned at the beginning of the paper.

Let $M$ be the lattice of $A_2.$ Let $L={\mathbb Z} \lambda \oplus
M^6.$ Then $L$ has a primitive embedding into the $K3$ lattice. By
Theorem \ref{urabe} there is a sextic curve $C_1$ such that the
Picard group of the corresponding double sextic is isomorphic to
$L.$

Denote the $12$ generators of the root lattice of $6A_2$ by $e_i (1
\leq i \leq 12)$ such that $e_i e_{i+1}=1$ for $i=1,3,5,7,9,11.$ Let
$$u=\sum_{i=1}^6 \frac{e_{2i-1}+2e_{2i}}{3}.$$ Then $u v \in {\mathbb Z}$
for every $v \in L$ and $u^2 \in 2{\mathbb Z}.$ Hence the subgroup
$L'$ of $L^\vee$ generated by $u$ and $L$ is an overlattice of $L.$
By Nikulin's embedding criterion there is a primitive embedding of
$L'$ into $\Lambda$. It is easy to check that the two additional
conditions of Theorem \ref{urabe} are also satisfied. Thus there is
a sextic curve $C_2$ such that the Picard group of the corresponding
double sextic is isomorphic to $L'.$

Both $C_1$ and $C_2$ are irreducible sextic curves with six $A_2$
cusps as their only singularities. However, their discriminantal
groups are $L^\vee/L$ and ${L'}^\vee/L'$ respectively, which are not
isomorphic. Hence $\{C_1,C_2\}$ is a Zariski pair.

Next we show that the six cusps of $C_2$ are located on a conic.

Let $p: Y \rightarrow {\mathbb P}^2$ be the double cover of
${\mathbb P}^2$ branched over $C_2$ and let $\mu: X \rightarrow Y$
be the minimal resolution of singularities of $Y.$ Identify $\text{
Pic}(X)$ with $L'.$ The map $p \mu$ is determined by the linear
system $|\lambda|.$  Let
$$D=\lambda-\sum_{i=1}^6 \frac{e_{2i-1}+2e_{2i}}{3}.$$
Then $D \in \text{Pic}(X).$ Since $D^2=-2,$ the Riemann-Roch theorem
implies that either $h^0(X,D)>0$ or $h^0(X,-D)>0.$ Since $\lambda D
=2 >0,$ we have $h^0(X,D)>0.$ Hence we may assume that $D$ is an
effective divisor. Choose an irreducible component $D_1$ of $D$ such
that $\lambda D_1 >0.$

Suppose that $\lambda D_1 =1.$ Then the divisor $D_1$  would be in
the class $\lambda/2 + \sum_{i=1}^{12} k_i e_i,$ in which $k_i \in
{\mathbb Q}.$ However, the latter is not in $\text{Pic}(X).$ This
leads to a contradiction. Hence $\lambda D_1 = 2.$

Let $D=D_1+E.$ Then $\lambda F=0$ for each irreducible component $F$
of $E.$ Thus $F$ is contracted to a point in ${\mathbb P}^2.$ This
means that $E=\sum_{j=1}^{12} k_j e_j$ where $k_j$ is a nonnegative
integer. For any $1 \leq i \leq 6,$ if $k_{2i-1} >0$ or $k_{2i}>0$
then
\begin{eqnarray*}
0 & \leq & k_{2i}
=  D(k_{2i-1}e_{2i-1}+k_{2i}e_{2i}) \\
& = &
D_1(k_{2i-1}e_{2i-1}+k_{2i}e_{2i})+(k_{2i-1}e_{2i-1}+k_{2i}e_{2i})^2
\\ & < & D_1(k_{2i-1}e_{2i-1}+k_{2i}e_{2i})
\end{eqnarray*}
implies that either $D_1e_{2i-1}>0$ or $D_1e_{2i}>0.$
 Hence $p \mu(D_1)$ is a conic passing through all six cusps
of $C_2.$

\section{Sextics of Milnor number $19$ with simple singularities}

Shimada finds all arithmetic Zariski pairs for sextics with simple
singularities of Milnor number $19$ (\cite{Shimada1}). In this
section we present a few more Zariski pairs of sextics of Milnor
number $19$ with different discriminantal groups, all reducible.

\medskip

{\bf Example 1} $E_6+A_{11}+2A_1$

\medskip

\begin{figure}
    \centering
    \begin{minipage}[t]{0.8\linewidth}
        \centering
        \includegraphics[scale=1]{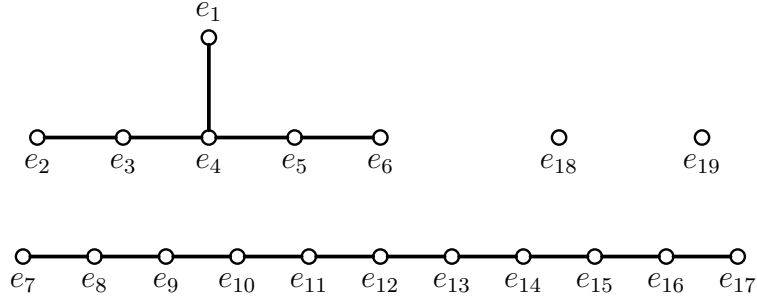}
        \caption{Dynkin graph of $E_6+A_{11}+2A_1$}
        \label{fig01}
    \end{minipage}
\end{figure}

Let $L$ denote the negative definite lattice of the Dynkin graph
$E_6+A_{11}+2A_1.$ The $19$ generators of $L$ are labeled according
to Figure \ref{fig01}. Let $V={\mathbb Q} \otimes_{\mathbb Z}(
{\mathbb Z}\lambda \oplus L),$ in which $\lambda^2=2.$ Let
$$u=\frac{\sum_{i=1}^{11}
ie_{i+6}}{2}+\frac{e_{18}}{2}+\frac{e_{19}}{2} \in V.$$

It can be verified that $u^2 \in 2{\mathbb Z}$ and $u w \in {\mathbb
Z}$ for any $w \in {\mathbb Z} \lambda \oplus L.$ Let $M_1$ be the
lattice generated by ${\mathbb Z}\lambda \oplus L$ and $u.$ Then
$M_1$ is an overlattice of ${\mathbb Z}\lambda \oplus L.$ Using
Nikulin's criterion for lattice embedding (\cite{nik}1.12.2), one
verifies that there is a primitive embedding from $M_1$ into the K3
lattice $\Lambda.$ Moreover, it is not hard to check that $M_1$
satisfies the two additional conditions in Theorem \ref{urabe}. It
follows that there is a reduced sextic curve $C_1$ with
$E_6,A_{11},A_1,A_1$ as its singularities. Although we can use the
algorithm in \cite{yang} to determine the irreducible decomposition
of $C_1,$ the following lemma uses an elementary argument to serve
the same purpose.

\begin{lemma}\label{irr01}
Let $C$ be a reduced sextic curve with $E_6,A_{11},A_1,A_1$ as its
only singularities. Then $C=B+D$ where $B$ and $D$ are irreducible
curves of degree $2$ and degree $4$ respectively satisfying the
following conditions:

1) $D$ has an $E_6$ singularity;

2) $B \cap D =\{p,q_1,q_2\},$ in which $p$ is an $A_{11}$ point of
$C$ and $B,D$ meet at $q_1$ and $q_2$ transversally.
\end{lemma}

Proof. Let $d$ be the maximal degree of all irreducible components
of $C.$ Since $E_6$ is locally irreducible, $d$ is at least $4.$
Since the arithmetic genus of an irreducible sextic curve is $10,$
there is no irreducible sextic curve with $E_6,A_{11},A_1,A_1$ as
its singularities. Hence $d \leq 5.$ It is obvious that there are
two irreducible components passing through the $A_{11}$ point and
the intersection number of these two components at $A_{11}$ point is
$6.$ The only possibility is that $d=4$ and the other component is a
conic. The rest is clear. $\Box$

\medskip

Let
$$v=\frac{3e_1+2e_2+4e_3+6e_4+5e_5+4e_6}{3}+\frac{\sum_{i=1}^{11}
ie_{i+6}}{6}+\frac{e_{18}}{2}+\frac{e_{19}}{2} \in V.$$

Then $v^2 \in 2{\mathbb Z}$ and $v w \in {\mathbb Z}$ for any $w \in
{\mathbb Z} \lambda \oplus L.$ Let $M_2$ be the lattice generated by
${\mathbb Z}\lambda \oplus L$ and $v.$ Then $M_2$ is an overlattice
of ${\mathbb Z}\lambda \oplus L.$ Note that $3v \equiv u \pmod{L}.$

Using the same method as before, we assert that $M_2$ satisfies the
conditions in Theorem \ref{urabe}. Let $C_2$ be the sextic
determined by $M_2$. By Lemma \ref{irr01} $C_2$ has the same
configuration as $C_1.$ The discriminantal groups of $C_1$ and $C_2$
are $M_1^\vee/M_1$ and $M_2^\vee/M_2$ respectively. They are finite
groups of different size. Hence they are not isomorphic. This shows
that $\{C_1,C_2\}$ is a Zariski pair.

\medskip

Next we show that $C_1$ and $C_2$ are distinguished by the existence
of a special conic on ${\mathbb P}^2.$

\begin{lemma}\label{conic01}
Let $C$ be either $C_1$ or $C_2.$ Let $p_1,p_2,p_3,p_4$ be its
singularities of types $E_6,A_{11},A_1,A_1$ respectively.  Then the
following statements hold:

1) There is a conic $Q$ on ${\mathbb P}^2$ such that
$(C_2,Q)_{p_1}=(C_2,Q)_{p_2}=4, (C_2,Q)_{p_3}=(C_2,Q)_{p_4}=2;$

2) There is no conic $Q$ on ${\mathbb P}^2$ such that
$(C_1,Q)_{p_1}=(C_1,Q)_{p_2}=4, (C_1,Q)_{p_3}=(C_1,Q)_{p_4}=2.$
\end{lemma}

Proof. Let $p: Y \rightarrow {\mathbb P}^2$ be the double cover
branched over $C$ and let $\mu: X \rightarrow Y$ be the minimal
resolution of singularities of $Y.$ There are nineteen $-2$ curves
on $X$ arising from the singularities of $C.$ They are still denoted
by $e_1,\ldots,e_{19}$ by abuse of notation.

\medskip

Case 1) $C=C_2:$

The map $p \mu$ is determined by the linear system $|\lambda|.$ Let
$$D=\lambda-\frac{3e_1+2e_2+4e_3+6e_4+5e_5+4e_6}{3}-\frac{\sum_{i=1}^{11}
ie_{i+6}}{6}+e_{17}-\frac{e_{18}}{2}-\frac{e_{19}}{2}.$$ It follows
from $D \equiv -v \pmod{{\mathbb Z}\lambda \oplus L}$ that $D \in
\text{Pic}(X).$ Since $D^2=-2,$ the Riemann-Roch theorem shows that
either $h^0(X,D)>0$ or $h^0(X,-D)>0.$ Since $\lambda D =2
>0,$ we have $h^0(X,D)>0.$ Hence we may assume that $D$ is an
effective divisor. Choose an irreducible component $D_1$ of $D$ such
that $\lambda D_1 >0.$

Suppose that $\lambda D_1 =1.$ Then image of $D_1$ in $L^\vee/L$
would be $\lambda/2+\sum_{i=1}^{19} n_i e_i$ for some $n_i(1 \leq i
\leq 19),$ which is not in the lattice $M_2.$ This leads to a
contradiction. Hence $\lambda D_1 = 2.$

Let $D=D_1+E.$ Then $\lambda F=0$ for each irreducible component $F$
of $E.$ Thus $F$ is contracted to a point in ${\mathbb P}^2.$ This
means that $E$ consists of exceptional curves. Since $D e_i >0$ for
$i=6,16,18,19,$ the image $Q$ of $D$ in ${\mathbb P}^2$ passes
through all singularities of $C_2.$ Since there is no line with this
property, it must be a conic.

Following the process of the canonical resolution of a double cover,
it is easy to see that $(C_2,Q)_{p_1}=(C_2,Q)_{p_2}=4.$

\medskip

Case 2) $C=C_1:$

Suppose that there is such a conic $Q$ for $C_1.$ Since $Q$ has even
intersection number with $C_1$ at each point of $Q \cap C_1,$ it
splits into two components $\tilde{Q}_1$ and $\tilde{Q}_2$ in $X.$
Since $\lambda \tilde{Q}_i>0$ for $i=1,2$ and $\lambda
(\tilde{Q}_1+\tilde{Q}_2)=2,$ we have $\lambda \tilde{Q}_1=\lambda
\tilde{Q}_2=1.$ Since $(Q,C_1)_{p_2}=4,$ one of $\tilde{Q}_1$ and
$\tilde{Q}_2$ meets $e_{16}$ transversally. We may assume that
$\tilde{Q}_1 e_{16}=1.$ Since $\tilde{Q}_1 e_i=0$ for $7 \leq i \leq
15$ and $i=18,$ the divisor class $[\tilde{Q}_1]$ is equal to
$$\lambda-\sum_{i=1}^6 n_i e_i
-(e_7+2e_8+3e_9+4e_{10}+5e_{11}+6e_{12}+7e_{13}+8e_{14}+9e_{15}+10e_{16}+5e_{17})/6-\sum_{i=18}^{19}
m_i e_i$$ for some $n_i$ and $m_i.$ Obviously there is no such
element in $M_1.$ $\Box$

\medskip

The equations are as follows:

$$C_1:
(3x_0^3x_2+3x_0^2x_1^2-3x_0x_1^3+2x_1^4)(3x_0x_2-x_2^2+3x_1x_2+3x_1^2)=0,$$

for which $E_6,A_{11},A_1,A_1$ are located at
$(0:0:1),(1:0:0),(1:3/2+\sqrt{-3}:33/4+15\sqrt{-3}/4),(1:3/2-\sqrt{-3}:33/4-15\sqrt{-3}/4)$
respectively.

$$C_2:
(x_1x_2+3x_0x_1+8x_0^2+3x_0x_2)(3x_1^3x_2+x_0x_1^3-2x_1^2x_2^2$$
$$+3x_0x_1^2x_2+3x_1x_2^3+3x_0x_1x_2^2+x_0x_2^3)=0,$$

for which $E_6,A_{11},A_1,A_1$ are located at
$(1:0:0),(1:-2:-2),(0:1:0),(0:0:1)$ respectively.

For $C_2$ the conic $Q: x_0x_1+x_0x_2+x_1x_2=0$ passes through all
singularities. The intersection numbers of $Q$ with the quartic at
$E_6$ and $A_{11}$ are equal to $4.$

\medskip

{\bf Example 2} $A_{11}+A_5+3A_1$

\medskip

\begin{figure}
    \centering
    \begin{minipage}[t]{0.8\linewidth}
        \centering
        \includegraphics[scale=1]{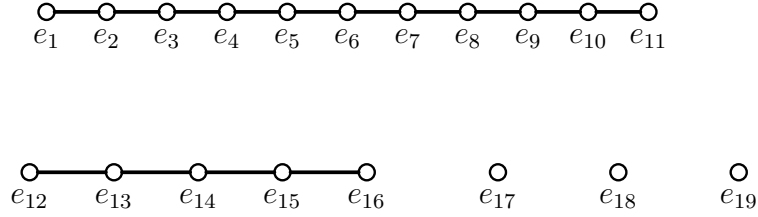}
        \caption{Dynkin graph of $A_{11}+A_5+3A_1$}
        \label{fig02}
    \end{minipage}
\end{figure}

Let $L$ denote the lattice of the Dynkin graph $A_{11}+A_5+3A_1.$
The $19$ generators of $L$ are labeled according to Figure
\ref{fig02}. Let $V={\mathbb Q} \otimes_{\mathbb Z}( {\mathbb
Z}\lambda \oplus L).$  Let
$$u_1=\frac{\sum_{i=1}^{11}
ie_i}{2}+\frac{e_{18}}{2}+\frac{e_{19}}{2},$$
$$u_2=\frac{\lambda}{2}+\frac{\sum_{i=1}^{5}
ie_{i+11}}{2}+\frac{e_{18}}{2}+\frac{e_{19}}{2}.$$ Let $M_1$ be the
lattice generated by $u_1,u_2$ and ${\mathbb Z}\lambda \oplus L.$
Then $M_1$ is an overlattice of ${\mathbb Z}\lambda \oplus L$
satisfying the conditions in Theorem \ref{urabe}. The configuration
of $C_1$ can be described as follows. $C_1$ consists of three
components $Q_1,Q_2,Q_3$ of degrees $1,2,3$ respectively. $Q_3$ has
a node. $Q_2 \cap Q_3$ is an $A_{11},$ $Q_1 \cap Q_3$ is an $A_3$
point and $Q_1$ intersects $Q_2$ at two distinct points.

Let
$$v=\frac{\sum_{i=1}^{11}
ie_i}{6}+\frac{\sum_{i=1}^{5}
ie_{i+11}}{3}+\frac{e_{18}}{2}+\frac{e_{19}}{2}.$$ Let $M_2$ be the
lattice generated by $u_1,u_2,v$ and ${\mathbb Z}\lambda \oplus L.$
Then $M_2$ is an overlattice of ${\mathbb Z}\lambda \oplus L$
satisfying the conditions in Theorem \ref{urabe}. The corresponding
sextic curves $C_2$ has the same configuration as $C_1.$ However
their discriminantal groups are not isomorphic, since their orders
are different. It follows that $\{C_1,C_2\}$ is a Zariski pair.

The equations are as follows:

$$C_1:
(2x_0-5x_1+3x_2)(16x_0x_2-16x_1^2-4x_1x_2+3x_2^2)(4x_0^2x_2-4x_0x_1^2-2x_0x_1x_2+x_1^3+x_1^2x_2)=0,$$

in which the $A_{11},A_5,A_1,A_1,A_1$ points are located at
$(1:0:0),(1:1:1),(0:0:1),(1:\frac{24}{43}+\frac{4
\sqrt{-3}}{43}:\frac{28}{43}+\frac{20
\sqrt{-3}}{129}),(1:\frac{24}{43}-\frac{4
\sqrt{-3}}{43}:\frac{28}{43}-\frac{20 \sqrt{-3}}{129})$
respectively.

$$C_2:
(x_2-x_0)(x_0^2x_2-x_0x_1^2+x_1^2x_2)(x_0x_2-x_1^2+x_2^2)=0,$$

in which the $A_{11},A_5,A_1,A_1,A_1$ points are located at
$(1:0:0),(0:1:0),(0:0:1),(1:\sqrt{2}:1),(1:-\sqrt{2}:1)$
respectively.

\medskip

For $C_2$ the line $L: x_2=0$ passes $A_{11},A_5$ and is tangent to
the conic at $A_{11}.$ This property is not shared by $C_1.$

\medskip

{\bf Example 3} $A_{17}+2A_1$

\medskip

\begin{figure}
    \centering
    \begin{minipage}[t]{0.8\linewidth}
        \centering
        \includegraphics[scale=1]{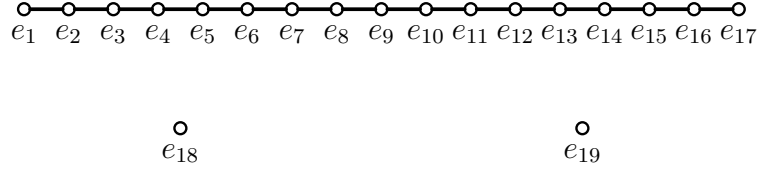}
        \caption{Dynkin graph of $A_{17}+2A_1$}
        \label{fig03}
    \end{minipage}
\end{figure}

Let $L$ denote the lattice of the Dynkin graph $A_{17}+2A_1.$ The
$19$ generators of $L$ are labeled according to Figure \ref{fig03}.
Let $V={\mathbb Q} \otimes_{\mathbb Z}( {\mathbb Z}\lambda \oplus
L).$  Let
$$u=\frac{\lambda}{2}+\frac{\sum_{i=1}^{17}
ie_i}{2}.$$

Let $M_1$ be the lattice generated by $u$ and ${\mathbb Z}\lambda
\oplus L.$ Then $M_1$ is an overlattice of ${\mathbb Z}\lambda
\oplus L$ satisfying the conditions in Theorem \ref{urabe}. Let
$C_1$ be the corresponding sextic.

Let
$$v=\frac{\lambda}{2}+\frac{\sum_{i=1}^{17}
ie_i}{6}.$$

Let $M_2$ be the lattice generated by $v$ and ${\mathbb Z}\lambda
\oplus L.$ Then $M_2$ is an overlattice of ${\mathbb Z}\lambda
\oplus L$ satisfying the conditions in Theorem \ref{urabe}. Let
$C_2$ be the corresponding sextic.

Both curves $C_1$ and $C_2$ have two nodal cubics as irreducible
components and these two components meet at $A_{17}.$ They form a
Zariski pair.

Denote the $A_{17}$ point by $p.$ Let $L$ be the common tangent line
of the two cubics at $p.$ Then $(L,C_1)_p=4$ and $(L,C_2)_p=6.$

\medskip

The equations are as follows:

$$C_1:
(x_0^2x_2-x_0x_1^2+x_1^3+x_0x_1x_2+7x_1^2x_2)
[(29+9\sqrt{-3})x_0^2x_2-(29+9\sqrt{-3})x_0x_1^2$$
$$-(484-18\sqrt{-3})x_0x_1x_2+542x_1^3+(3078-54\sqrt{-3})x_0x_2^2
-(3901-135\sqrt{-3})x_1^2x_2$$
$$+(13851-243\sqrt{-3})x_2^3+(1539-27\sqrt{-3})x_1x_2^2]=0,$$

in which the $A_{17},A_1,A_1$ are located at $(1:0:0),(0:0:1)$ and
$(1:\frac{1}{13}+\frac{3\sqrt{-3}}{65}:\frac{1}{390}+\frac{7\sqrt{-3}}{1170})$
respectively.

$$C_2:
(x_0^2x_2+x_1^3+x_0x_1x_2+x_1^2x_2-\frac{x_2^3}{16})(x_0^2x_2+x_1^3+x_0x_1x_2+x_1^2x_2)=0,$$

in which the $A_{17},A_1,A_1$ are located at
$(1:0:0),(1:-2:4),(0:0:1)$ respectively.

\medskip

{\bf Example 4} $A_{15}+A_3+A_1$

\medskip

\begin{figure}
    \centering
    \begin{minipage}[t]{0.8\linewidth}
        \centering
        \includegraphics[scale=1]{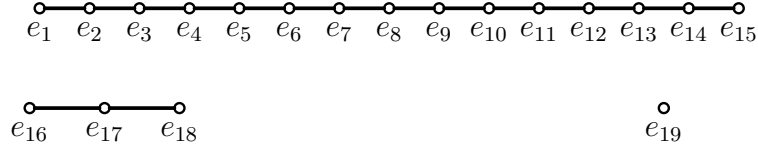}
        \caption{Dynkin graph of $A_{15}+A_3+A_1$}
        \label{fig04}
    \end{minipage}
\end{figure}

Let $L$ denote the lattice of the Dynkin graph $A_{15}+A_3+A_1.$ The
$19$ generators of $L$ are labeled according to Figure \ref{fig04}.
Let $V={\mathbb Q} \otimes_{\mathbb Z}( {\mathbb Z}\lambda \oplus
L).$  Let
$$u=\frac{\lambda}{2}+\frac{\sum_{i=1}^{15}
ie_i}{4}+\frac{e_{16}+2e_{17}+3e_{18}}{2}+\frac{e_{19}}{2}.$$

Let $M_1$ be the lattice generated by $u$ and ${\mathbb Z}\lambda
\oplus L.$ Then $M_1$ is an overlattice of ${\mathbb Z}\lambda
\oplus L$ satisfying the conditions in Theorem \ref{urabe}. Let
$C_1$ be the corresponding sextic.

Let
$$v=\frac{\lambda}{2}+\frac{\sum_{i=1}^{15}
ie_i}{8}+\frac{e_{16}+2e_{17}+3e_{18}}{4}.$$

Let $M_2$ be the lattice generated by $u,v$ and ${\mathbb Z}\lambda
\oplus L.$ Then $M_2$ is an overlattice of ${\mathbb Z}\lambda
\oplus L$ satisfying the conditions in Theorem \ref{urabe}. Let
$C_2$ be the corresponding sextic.

Both curves $C_1$ and $C_2$ have two irreducible components of
degree $4$ and $2$ respectively. The quartic component contains an
$A_3$ and an $A_1.$ The two components meet at $A_{15}.$
$\{C_1,C_2\}$ is a Zariski pair.

\medskip

The equations are as follows:

$$C_1:
(2x_0x_2+2x_1^2-2x_2^2+\sqrt{-2}x_1x_2)(4x_0^3x_2+4x_0^2x_1^2+$$
$$4x_0^2x_2^2-2\sqrt{-2}x_0^2x_1x_2-4\sqrt{-2}x_0x_1^3
+34x_0x_1^2x_2+$$
$$12\sqrt{-2}x_0x_1x_2^2+22x_1^4+15\sqrt{-2}x_1^3x_2-18x_1^2x_2^2)=0,
$$

in which $A_{15},A_3$ and $A_1$ are located at $(1:0:0),(0:0:1)$ and
$(1:-\frac{\sqrt{-2}}{3}:\frac{2}{9})$ respectively. The line
$x_1=0$ passes $A_{15}$ and $A_3$ but its intersection number with
$C_1$ at $A_{15}$ is equal to $2.$

$$C_2:
(x_0x_2+x_1^2-\frac{x_2^2}{2})(x_0^3x_2+x_0^2x_1^2+2x_0x_1^2x_2+x_1^2x_2^2+\frac{3x_0^2x_2^2}{2})=0,$$

in which $A_{15},A_3$ and $A_1$ are located at $(1:0:0),(0:1:0)$ and
$(0:0:1)$ respectively. The intersection numbers of the line $x_2=0$
with $C_2$ at $A_{15}$ and $A_3$ are $4$ and $2$ respectively.

\medskip

{\bf Example 5} $2A_9+A_1$

\medskip

\begin{figure}
    \centering
    \begin{minipage}[t]{0.8\linewidth}
        \centering
        \includegraphics[scale=1]{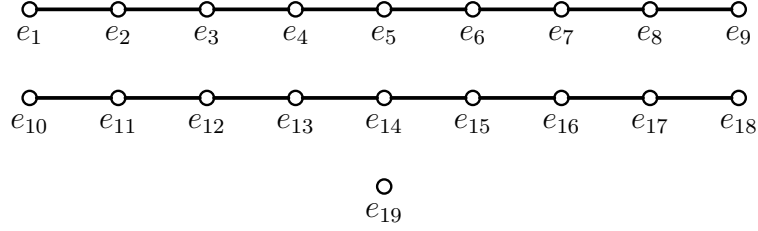}
        \caption{Dynkin graph of $2A_9+A_1$}
        \label{fig05}
    \end{minipage}
\end{figure}

Let $L$ denote the lattice of the Dynkin graph $2A_9+A_1.$ The
generators of $L$ are labeled according to Figure \ref{fig05}. Let
$V={\mathbb Q} \otimes_{\mathbb Z}( {\mathbb Z}\lambda \oplus L).$
Let
$$u=\frac{\lambda}{2}+\frac{\sum_{i=1}^{9}
ie_i}{2}.$$

Let $M_1$ be the lattice generated by $u$ and ${\mathbb Z}\lambda
\oplus L.$ Then $M_1$ is an overlattice of ${\mathbb Z}\lambda
\oplus L$ satisfying the conditions in Theorem \ref{urabe}. Let
$C_1$ be the corresponding sextic.

Let
$$v=\frac{\lambda}{2}+\frac{\sum_{i=1}^{9}
ie_i}{5}+\frac{\sum_{i=1}^{9} ie_{i+9}}{10}.$$

Let $M_2$ be the lattice generated by $u,v$ and ${\mathbb Z}\lambda
\oplus L.$ Then $M_2$ is an overlattice of ${\mathbb Z}\lambda
\oplus L$ satisfying the conditions in Theorem \ref{urabe}. Let
$C_2$ be the corresponding sextic.

Both curves $C_1$ and $C_2$ have two irreducible components of
degree $5$ and $1$ respectively. The quintic component contains an
$A_9$ and an $A_1.$ The two components meet at the second $A_9.$ Let
$L$ be the line connecting the two $A_9$ points. Let $p$ be the
$A_9$ of the quintic component. Then $(L,C_1)_p=2$ and
$(L,C_2)_p=4.$

\medskip

The equations of $C_1$ and $C_2$ are as follows:

$$C_1:
x_0[x_0^3x_2^2+2x_0^2x_1^2x_2+x_0x_1^4+(25/2-5\sqrt{5}/2)x_1^5+(27-5\sqrt{5})x_0x_1^3x_2$$
$$+(-22+6\sqrt{5})x_0^2x_2^3+(29/2-5\sqrt{5}/2)x_0^2x_1x_2^2+(-21+6\sqrt{5})x_0x_1^2x_2^2$$
$$+(-19/2+7\sqrt{5}/2)x_0x_1x_2^3+(-9+4\sqrt{5})x_0x_2^4]=0,
$$

in which the quint component has an $A_9$ at $(1:0:0)$ and an $A_1$
at $(1:-39-87\sqrt{5}/5:648+1449\sqrt{5}/5),$ and the two components
meet at another $A_9$ point $(0:1:0).$

$$C_2:
x_1(x_0^2x_1^3+2x_0x_1^2x_2^2+x_1x_2^4+\frac{3x_0^3x_1^2}{4}
+x_0^2x_1^2x_2+\frac{9x_0^4x_1}{64}$$
$$+\frac{3x_0^3x_1x_2}{8}+x_0^2x_1x_2^2+x_0x_1x_2^3+\frac{125x_0^5}{27648})=0,
$$

in which $A_1$ is at $(1:-5/48:-1/4)$ and the two $A_9$ points are
at $(0:1:0)$ and $(0:0:1)$ respectively.

In summary, we obtain the following result.

\begin{thm}\label{maximizing}
There are five Zariski pairs of sextic curves of Milnor number $19,$
whose combinations of singularities are
$$E_6+A_{11}+2A_1, A_{17}+2A_1, A_{15}+A_3+A_1, A_{11}+A_5+3A_1,
2A_9+A_1.$$
\end{thm}

\section{Sextic Zariski pairs with Milnor number less than $19$}

Many Zariski pairs of sextics with lower Milnor numbers can be
obtained using the same method. We choose several special ones to
analyze in details. The others can be obtained by the same method
and are listed in the tables at the end of this paper.

\subsection{Zariski triplets}

Four Zariski triplets were found among reduced sextics with simple
singularities.

\subsubsection{Zariski triplet of three conics}

Namba and Tsuchihashi constructed a Zariski pair of octive curves
which consists of four conics in \cite{Namba}. Here we give a
Zariski triplet consisting of three conics.

Let $L$ be the root lattice of $3A_5+3A_1.$ Let $N={\mathbb
Z}\lambda \oplus L.$ Denote the $18$ generators of $L$ by $e_i (1
\leq i \leq 18)$ in a natural way such that $e_ie_{i+1}=1$ for $1
\leq i \leq 4, 6 \leq i \leq 9$ and $11 \leq i \leq 14.$ Let
$$u_1=(\sum_{i=1}^5 i(e_{i+5}+e_{i+10})+e_{17}+e_{18})/2,$$
$$v_1=(\sum_{i=1}^5 i(e_i+e_{i+10})+e_{16}+e_{18})/2,$$
$$u_2=\sum_{i=1}^5 i(e_i+e_{i+5})/6+\sum_{i=1}^5 ie_{i+10}/3+(e_{17}+e_{18})/2,$$
$$v_2=(\sum_{i=1}^5 i(e_{i+5}+e_{i+10})+e_{16}+e_{18})/2$$
and
$$w=\lambda/2+\sum_{i=1}^5 i(e_i+e_{i+5}+e_{i+10})/2.$$

Let $A_i$ be the sublattice of $N^\vee$ generated by $u_i,v_i$ and
$N$ for $i=1,2.$ Let $A_3$ be the sublattice of $N^\vee$ generated
by $u_1,v_1$ and $w.$ Then $A_1,A_2,A_3$ are overlattices of $N$
satisfying the conditions in Theorem \ref{urabe}. Let $E_1,E_2,E_3$
be their corresponding reduced sextic curves. The configurations of
these two curves turn out to be the same: three conics. Since the
orders of $A_1^\vee/A_1,A_2^\vee/A_2$ and $A_3^\vee/A_3$ are
different, $\{E_1,E_2,E_3\}$ is a Zariski triplet.

Let us call a sextic curve corresponding to the overlattice $A_i$ a
sextic of type $i$ for $i=1,2,3.$

\begin{thm}\label{threeconics}
There is a Zariski triplet $\{E_1,E_2,E_3\}$ of sextics of three
conics where $E_i$ is of type $i.$ They are distinguished by the
following conditions:

1) For $E_2$ there is a conic on ${\mathbb P}^2$ passing through the
three $A_5$ points such that the intersection number of this conic
with $E_2$ at each $A_5$ point is $4.$ This property is not shared
by $E_1$ and $E_3.$

2) For $E_3$ there is a nodal cubic with the node at one $A_5$ point
such that the intersection number of this nodal cubic with $E_3$ at
each $A_5$ point is equal to $6.$ This property is not shared by
$E_1$ and $E_2.$
\end{thm}

Sketch of the proof. Let $X_i$ denote the $K3$ surface obtained by
the double cover branched over $E_i$ for $i=1,2,3.$

Let
$$D=\lambda-\frac{e_1+2e_2+3e_3+4e_4+2e_5}{3}
-\frac{e_6+2e_7+3e_8+4e_9+2e_{10}}{3}$$
$$-\frac{2e_{11}+4e_{12}+3e_{13}+2e_{14}+e_{15}}{3},$$ Then $D \in
\text{Pic}(X_2)$ and $D^2=-2, D \lambda = 2.$ Using the same
argument as before, we proved that the image of one irreducible
component of a member of $|D|$ under the double cover map is a conic
with the desired property.

Let
\begin{eqnarray*}
C & = &
3\lambda/2-\frac{e_1+2e_2+3e_3+2e_4+e_5}{2}-\frac{e_6+2e_7+3e_8+2e_9+e_{10}}{2}\\
& & -\frac{e_{11}+2e_{12}+3e_{13}+4e_{14}+3e_{15}}{2},
\end{eqnarray*}
 Then $C \in
\text{Pic}(X_3)$ and $C^2=-2, C \lambda =3.$ For the same reason as
before the image of one irreducible component of a member of $|C|$
in ${\mathbb P}^2$ is a nodal cubic with the desired property.
$\Box$

\medskip

In the remaining part of this subsection we calculate in details the
explicit equations of all such sextics.

Let $(x_0:x_1:x_2)$ be the homogeneous coordinates of ${\mathbb
P}^2.$ Let $Y$ be a sextic curve consisting of three conics
$C_1,C_2,C_3$ satisfying the conditions in Theorem
\ref{threeconics}. After a suitable linear change of coordinates, we
may assume that the three $A_5$ points are located at
$$(1:0:0) \in C_2 \cap C_3,(0:1:0) \in C_1 \cap C_3,(0:0:1) \in C_1 \cap C_2$$
and the tangent lines of $C_1$ at $(0:1:0)$ and $(0:0:1)$ are given
by the equations $x_2-x_0=0$ and $x_1-x_0=0$ respectively. The
equations of the conics are written as

$$C_1: x_1x_2-x_0x_1-x_0x_2+\lambda x_0^2=0,$$
$$C_2: x_1x_2-ax_0x_1-x_0x_2+bx_1^2=0,$$
$$C_3: x_1x_2-x_0x_1-\frac{1}{a}x_0x_2+cx_2^2=0$$

where $\lambda,a,b,c$ are parameters to be determined. With generic
values for $\lambda, a,b,c$ the sextic has at least three $A_3$
singularities at $(1:0:0),(0:1:0)$ and $(0:0:1)$ already. The
requirement of $A_5$ poses three conditions on the parameters, which
are determined as follows.

1) $A_5$ at $(1:0:0)$:

Under the standard affine coordinates $x=x_1/x_0,y=x_2/x_0$ the
equations of $C_2$ and $C_3$ are

$$C_2: xy-ax-y+bx^2=0,$$
$$C_3: xy-x-\frac{1}{a}y+cy^2=0.$$

After the change of coordinates $y=y'-ax$ these two equations become

\begin{equation}\label{eq1}
y'+(a-b)x^2-xy'=0
\end{equation}

and

\begin{equation}\label{eq2}
y'+a^2(1-ac)x^2-a(1-2ac)xy'+cy'^2=0.
\end{equation}

In order that the intersection number of $C_2$ and $C_3$ at
$(1:0:0)$ is greater than two, the coefficients of the term $x^2$ in
(\ref{eq1}) and (\ref{eq2}) should be equal. This gives the
following relation

\begin{equation}\label{eq3}
b=a-a^2+a^3c.
\end{equation}

2) $A_5$ at $(0:1:0)$:

Under the affine coordinates $w=x_0/x_1,y=x_2/x_1$ the equations of
$C_1$ and $C_3$ are

$$C_1: y-w-wy+\lambda w^2=0,$$
$$C_3: y-w-\frac{1}{a}wy+cy^2.$$

After the change of coordinates $y=y'+w$ these two equations become

\begin{equation}\label{eq4}
y'+(\lambda-1)w^2-wy'=0
\end{equation}

and

\begin{equation}\label{eq5}
y'+(c-\frac{1}{a})w^2+\cdots =0.
\end{equation}

The coefficients of the term $w^2$ in (\ref{eq4}) and (\ref{eq5})
should be equal. This gives the following relation

\begin{equation}\label{eq6}
c=\frac{1}{a}+\lambda-1.
\end{equation}

3) $A_5$ at $(0:0:1)$:

Under the affine coordinates $w=x_0/x_2,x=x_1/x_2$ the equations of
$C_1$ and $C_3$ are

$$C_1: x-xw-w+\lambda w^2=0,$$
$$C_2: x-axw-w+bx^2=0.$$

The requirement of an $A_5$ at $(0:0:1)$ gives the following
condition

\begin{equation}\label{eq7}
b=a+\lambda-1.
\end{equation}

The three conditions (\ref{eq3}),(\ref{eq6}) and (\ref{eq7}) yield
$$(\lambda-1)(a^3-1)=0.$$
Since $\lambda \neq 1$ (otherwise the conic $C_1$ would become the
union of two lines), $a$ is equal to one of $1,\zeta,\zeta^2$ with
$\zeta=e^{2\pi i/3}.$ Thus we obtain three sets of solutions for the
parameters $a,b,c:$

$$a=1,b=\lambda,c=\lambda,$$
$$a=\zeta,b=\lambda+\zeta-1,c=\lambda+\zeta^2-1,$$
$$a=\zeta^2,b=\lambda+\zeta^2-1,c=\lambda+\zeta-1,$$

Hence there are three families of sextic curves satisfying our
conditions defined by

$$(x_1x_2-x_0x_1-x_0x_2+\lambda x_0^2)(x_1x_2-x_0x_1-x_0x_2+\lambda
x_1^2)$$

\begin{equation}\label{family1}
(x_1x_2-x_0x_1-x_0x_2+\lambda x_2^2)=0,
\end{equation}

$$(x_1x_2-x_0x_1-x_0x_2+\lambda x_0^2)(x_1x_2-\zeta x_0x_1-x_0x_2+ (\lambda
+\zeta-1) x_1^2)$$

\begin{equation}\label{family2}
(x_1x_2-x_0x_1-\zeta^2 x_0x_2+(\lambda+\zeta^2-1) x_2^2)=0,
\end{equation}

and

$$(x_1x_2-x_0x_1-x_0x_2+\lambda x_0^2)(x_1x_2-\zeta^2 x_0x_1-x_0x_2+
(\lambda +\zeta^2-1) x_1^2)$$

\begin{equation}\label{family3}
(x_1x_2-x_0x_1-\zeta x_0x_2+(\lambda+\zeta-1) x_2^2)=0
\end{equation}

respectively. The equation (\ref{family2}) becomes (\ref{family3})
if the variables $x_1$ and $x_2$ are exchanged. Hence they are
essentially the same.

Let $X_\lambda$ and $Y_\lambda$ denote the sextic curves defined by
(\ref{family1}) and (\ref{family2}) respectively. When $\lambda \neq
0,1,$ the sextic $X_\lambda$ is composed of three conics with
$3A_5+3A_1$ as singularities. The three $A_1$ points are located at
$$(1: \frac{2}{\lambda+1} : \frac{2}{\lambda+1}),(1: -1:
\frac{\lambda+1}{2}),(1: \frac{\lambda+1}{2}: -1).$$ When
$\lambda=-1,$ $X_\lambda$ is special in the sense that every two
$A_5$ points are collinear with an $A_1$ point.

When
$$\lambda \neq 0,\frac{3}{2},2,\frac{3 \pm \sqrt{-3}}{2},\frac{3 \pm \sqrt{-3}}{4},\frac{1 \pm
\sqrt{-3}}{2}$$ the sextic $Y_\lambda$ has the desired
configuration. The three $A_1$ points of $Y_\lambda$ are located at
$$\left(1: \frac{2\lambda}{2\lambda-3+\sqrt{-3}}:
\frac{\lambda(2\lambda-1+\sqrt{-3})}{4\lambda-3+\sqrt{-3}}\right),
\left(1: \frac{\lambda(2\lambda-1-\sqrt{-3})}{4\lambda-3-\sqrt{-3}}:
\frac{2\lambda}{2\lambda-3-\sqrt{-3}}\right)$$

and

$$\left(1:\frac{(-1+\sqrt{-3})(2\lambda-3)}{(\lambda-2)(2\lambda-3+\sqrt{-3})}:
\frac{(1-\sqrt{-3})(2\lambda-3)}{(\lambda-2)(2\lambda-3-\sqrt{-3})}\right).$$

\begin{thm}\label{special}
Let $X_\lambda$ and $Y_\lambda$ denote the sextic curves defined by
(\ref{family1}) and (\ref{family2}) respectively. Then

$$X_\lambda \left\{ \begin{array}{ll} \text{ is of type } 2, &
\text{ if } \lambda \neq 0,1 \\ \text{ degenerates, } & \text{ if }
\lambda = 0,1 \end{array} \right.$$

and

$$Y_\lambda \left\{ \begin{array}{ll} \text{ degenerates, } & \text{
if } \lambda = 0,\frac{3}{2},2,\frac{3 \pm \sqrt{-3}}{2},\frac{3 \pm
\sqrt{-3}}{4},\frac{1 \pm \sqrt{-3}}{2}, \\ \text{ is of type } 3, &
\text{ if } \lambda=3, \pm \sqrt{-3},\\ \text{ is of type } 1, &
otherwise. \end{array} \right.$$

In particular, every sextic of type 3 with this configuration is a
deformation of the ones of type 1.
\end{thm}

Proof.

The values of $\lambda$ for which $X_\lambda$ or $Y_\lambda$
degenerates have been discussed. Here we only consider the
non-degenerate values of $\lambda.$

Let $Q$ be the conic defined by $x_1x_2-x_0x_1-x_0x_2=0.$  It is
obvious that for every $X_\lambda$ or $Y_\lambda$ the conic $Q$ is
the unique one passing three $A_5$ points and whose intersection
numbers with the sextic at both $(0:1:0)$ and $(0:0:1)$ are greater
than $2.$ However, the intersection number of $Q$ with the sextic at
$(1:0:0)$ is greater than $2$ if and only if the sextic is some
$X_\lambda.$ Hence $X_\lambda$ is of type $2$ and $Y_\lambda$ is of
type $1$ or $3.$

Let $N$ be a nodal cubic characterizing some $Y_\lambda$ as a sextic
of type $3.$ There are three choices $(1:0:0),(0:1:0),(0:0:1)$ for
the location of the node. They will give three values of $\lambda.$
First assume that the node of $N$ is at $(1:0:0).$

Since the tangent lines of $N$ at $(0:1:0)$ and $(0:0:1)$ are the
same as those of $C_2,$ the equation of $N$ takes the form

$$x^2-x^2y+ry^2-rxy^2+sxy=0$$

under the affine coordinates $x=x_1/x_0,y=x_2/x_0.$ The tangent cone
of $N$ at $x=0,y=0$ is $x^2+ry^2+sxy$ and the tangent line of
$Y_\lambda$ is $y+\zeta x.$ In order that $N$ and $Y_\lambda$ have
higher contact at $x=0,y=0$ it is necessary that

\begin{equation}\label{eq8}
1-\zeta s + r \zeta^2 =0.
\end{equation}

The other two conditions of $N$ at $(0:1:0)$ and $(0:0:1)$ are
computed. They yields the following two relations:

\begin{equation}\label{eq9}
r-s=\lambda-1,
\end{equation}

\begin{equation}\label{eq10}
r(\lambda-1)-1+s=0.
\end{equation}

The solution of the simultaneous equations (\ref{eq8}),(\ref{eq9})
and (\ref{eq10}) is $r=1,s=-1,\lambda=3.$ Hence $Y_3$ has type $3.$

If the node is chosen to be at $(0:1:0)$ or $(0:0:1)$ then the
similar computation shows that $\lambda=\pm \sqrt{-3}.$ Therefore
$Y_3,Y_{\sqrt{-3}},Y_{-\sqrt{-3}}$ are all possible $Y_\lambda$ of
type $3.$  $\Box$

\subsubsection{$2A_7+A_3$}

Let $L$ be the root lattice of the Dynkin graph $2A_7+A_3.$ Let
$N={\mathbb Z}\lambda \oplus L.$ Denote the $17$ generators of $L$
by $e_i (1 \leq i \leq 17)$ in a natural way such that
$e_ie_{i+1}=1$ for $1 \leq i \leq 6, 8 \leq i \leq 13$ and $15 \leq
i \leq 16.$ Let
$$u_1=\sum_{i=1}^7 i(e_i+e_{i+7})/2,$$
$$u_2=\sum_{i=1}^7 i(e_i+e_{i+7})/4+(e_{15}+2e_{16}+3e_{17})/2,$$
$$u_3=\lambda/2+\sum_{i=1}^7 i(e_i+e_{i+7})/8+(e_{15}+2e_{16}+3e_{17})/4$$
be three elements in $N^\vee.$ It is easy to verify that $u_i^2 \in
2{\mathbb Z}$ for $i=1,2,3.$
 Let $M_i$ be the sublattice of
$N^\vee$ generated by $u_i$ and $N$ for $i=1,2,3.$ Note that $N
\subset M_1 \subset M_2 \subset M_3.$ It can be verified that all
these three overlattices of $N$ satisfy the conditions in Theorem
\ref{urabe}. Therefore there are three reduced sextic curves
$C_1,C_2,C_3$ whose discriminantal groups are
$M_1^\vee/M_1,M_2^\vee/M_2,M_3^\vee/M_3$ respectively. By the
algorithm in \cite{yang} we determine that all these three curves
have the same configuration. That is the union of an irreducible
quartic curve with an $A_3$ singularity and a conic such that the
conic meets the quartic at two points to form $2A_7.$ Thus we have
the following theorem

\begin{thm}\label{triple}
The collection of three sextic curves $\{C_1,C_2,C_3\}$ forms a
Zariski triplet.
\end{thm}

These three curves are distinguished by the following two
conditions:

1) The three singularities on $C_3$ are collinear while those on
$C_1$ or $C_2$ are not.

2) For $C_1$ and $C_2$ let $Q$ be the conic passing through the
three singularities and infinitely near points at both $A_7.$ Then
$Q$ passes the infinitely near point of $A_3$ for $C_2$ but does not
pass the one for $C_1.$

The equations are given as following:

\medskip

$$C_1:[x_0^2+ax_0x_1+(\sqrt{2}-a)x_0x_2+x_1x_2]
[2(\frac{\sqrt{2}}{2}a-1)^2x_0^2(x_2-x_1)^2+c_{2,2} x_1^2x_2^2$$
$$+c_{3,0} x_0x_1^3+c_{2,1}x_0x_1^2x_2
+c_{1,2}x_0x_1x_2^2-c_{0,3}x_0x_2^3+c_{3,1}x_1^3x_2+c_{1,3}x_1x_2^3]=0$$

where
$$c_{2,2}=2\sqrt{2} ab-2a^2b-b+2\sqrt{2}a -3,$$
$$c_{3,0}=a^3b-2\sqrt{2}a^2+2a,$$
$$c_{2,1}=3\sqrt{2}a^2b-3a^3b-2ab+4\sqrt{2}a^2-8a+2\sqrt{2},$$
$$c_{1,2}=3a^3b-6\sqrt{2}a^2b+8ab-2\sqrt{2}a^2-2\sqrt{2}b+6a-2\sqrt{2},$$
$$c_{0,3}=3\sqrt{2}a^2b-a^3b-6ab+2\sqrt{2}b,$$
$$c_{1,3}=a^2b-2\sqrt{2}ab+2b,$$
$$c_{3,1}=a^2b-2\sqrt{2}a+2,$$

in which $a$ and $b$ are generic parameters.

The $A_7,A_7,A_3$ are located at $(0:1:0),(0:0:1)$ and $(1:0:0)$
respectively.

\medskip

$$ C_2:(x_0^2+a x_0x_1-a x_0x_2+x_1x_2)
[a^2x_0^2(x_2-x_1)^2+(1-b-2a^2b)x_1^2x_2^2+a^3 b x_0x_1^3$$
$$+a(2-2b-3a^2b) (x_0x_1^2x_2-x_0x_1x_2^2)-a^3bx_0x_2^3+a^2bx_1^3x_2+a^2bx_1x_2^3 ]=0$$

for generic parameters $a$ and $b.$

The $A_7,A_7,A_3$ are located at $(0:1:0),(0:0:1)$ and $(1:0:0)$
respectively. The conic $ax_0x_1-ax_0x_2+x_1x_2$ is tangent to the
quartic component at $(1:0:0)$ and to the conic at both $(0:1:0)$
and $(0:0:1).$

\medskip

$$C_3:(x_0x_1+x_2^2)(x_0^3x_1+x_0^2x_2^2+bx_0^2x_1x_2-2x_0^2x_1^2+bx_0x_2^3$$
$$+\left(\frac{b^2}{4}-2\right)x_0x_1x_2^2-bx_0x_1^2x_2+x_0x_1^3+ax_2^4-bx_1x_2^3+x_1^2x_2^2)=0$$

for generic paramters $a$ and $b.$

The three collinear $A_7,A_7,A_3$ points are located
$(1:0:0),(0:1:0),(1:0:1)$ respectively.

\subsubsection{$2A_7+A_3+A_1$}

Let $L$ be the root lattice of the Dynkin graph $2A_7+A_3+A_1.$ Let
$N={\mathbb Z}\lambda \oplus L.$ Denote the $18$ generators of $L$
by $e_i (1 \leq i \leq 18)$ in a natural way such that
$e_ie_{i+1}=1$ for $1 \leq i \leq 6, 8 \leq i \leq 13$ and $15 \leq
i \leq 16.$ Let
$$u_1=\sum_{i=1}^7 i(e_i+e_{i+7})/2,$$
$$u_2=\sum_{i=1}^7 i(e_i+e_{i+7})/4+(e_{15}+2e_{16}+3e_{17})/2,$$
$$u_3=\sum_{i=1}^7 i(e_i+3e_{i+7})/8+(e_{15}+2e_{16}+3e_{17})/4+e_{18}/2$$
and
$$v=\lambda/2+\sum_{i=1}^7 ie_{i+7})/2+e_{18}/2$$
 be four elements in $N^\vee.$ It is easy to verify that $u_i^2
\in 2{\mathbb Z}$ for $i=1,2,3$ and $v^2 \in 2{\mathbb Z}.$
 Let $M_i$ be the sublattice of
$N^\vee$ generated by $u_i,v$ and $N$ for $i=1,2,3.$ It can be
verified that all these three overlattices of $N$ satisfy the
conditions in Theorem \ref{urabe}. Therefore there are three reduced
sextic curves $C'_1,C'_2,C'_3$ whose discriminantal groups are
$M_1^\vee/M_1,M_2^\vee/M_2,M_3^\vee/M_3$ respectively. By the
algorithm in \cite{yang} we determine that all these three curves
have the same configuration: the union of an irreducible quartic
curve with an $A_3$ singularity and two lines such that the each
line meets the quartic with intersection number $4.$ Thus we have
the following

\begin{thm}\label{triple}
The collection of three sextic curves $\{C'_1,C'_2,C'_3\}$ forms a
Zariski triplet.
\end{thm}

These three curves are distinguished by the following two
conditions:

1) The three singularities on $C'_3$ are collinear while those on
$C'_1$ or $C'_2$ are not.

2) For $C'_1$ and $C'_2$ let $Q$ be the conic passing through the
$A_7,A_7,A_3$ and infinitely near points at both $A_7.$ Then $Q$
passes the infinitely near point of $A_3$ for $C'_2$ but does not
pass the one for $C'_1.$

The equations of the curves are given as following:

\medskip

$$C'_1:(x_1-x_0)(x_2-x_0)
[\sqrt{-1}x_0^2x_1^2+2x_0^2x_1x_2-\sqrt{-1}x_0^2x_2^2-(2\lambda-2\sqrt{-1})x_0x_2^3$$
$$-(2\lambda+2-2\sqrt{-1})x_0x_1x_2^2-(2\lambda+2)x_0x_1^2x_2-2\lambda x_0x_1^3$$
$$+(2\lambda-2\sqrt{-1})x_1x_2^3+(2\lambda+2-\sqrt{-1})x_1^2x_2^2+2\lambda
x_1^3x_2]=0$$
for a generic parameter $\lambda.$

The $A_7,A_7,A_3$ and $A_1$ points are located at
$(0:1:0),(0:0:1),(1:0:0)$ and $(1:1:1)$ respectively.

\medskip

$$C'_2:(x_1-x_0)(x_2-x_0)
[x_0^2x_2^2+2x_0^2x_1x_2+x_0^2x_1^2-2\lambda
x_0x_1^3-(2\lambda+2)x_0x_1^2x_2$$
$$-(2\lambda+2)x_0x_1x_2^2-2\lambda x_0x_2^3+2\lambda x_1^3x_2+(2\lambda+1)x_1^2x_2^2+2\lambda x_1x_2^3]=0$$
for a generic parameter $\lambda.$

The $A_7,A_7,A_3,A_1$ points are located at
$(0:1:0),(0:0:1),(1:0:0)$ and $(1:1:1)$ respectively. The conic
$x_1x_2-x_0x_1-x_0x_2$ has intersection number $4$ with $C'_2$ at
$A_7$ and $A_3$ points.

\medskip

$$C'_3:x_1x_2(x_1^3x_2-2x_1^2x_2^2+x_1x_2^3+\lambda x_0^4)=0$$
for a generic parameter $\lambda.$

The $A_7,A_7,A_3,A_1$ points are located at
$(0:1:0),(0:0:1),(0:1:1)$ and $(1:0:0)$ respectively. The line
$x_0=0$ passes through the two $A_7$ points and the $A_3$ point.

\subsubsection{$A_{15}+A_3$}

Let $L$ be the root lattice of  $A_{15}+A_3.$ Let $N={\mathbb
Z}\lambda \oplus L.$ Denote the $18$ generators of $L$ by $e_i (1
\leq i \leq 18)$ in a natural way such that $e_ie_{i+1}=1$ for $1
\leq i  \leq 14$ and $16 \leq i \leq 17.$ Let
$$u_1=\sum_{i=1}^{15} \frac{ie_i}{2},$$
$$u_2=\sum_{i=1}^{15} \frac{ie_i}{4}+\frac{e_{16}+2e_{17}+3e_{18}}{2},$$
$$u_3=\lambda/2+\sum_{i=1}^{15} \frac{ie_i}{8}+\frac{e_{16}+2e_{17}+3e_{18}}{4}$$
be three elements in $N^\vee.$ It is easy to verify that $u_i^2 \in
2{\mathbb Z}$ for $i=1,2,3.$
 Let $M_i$ be the sublattice of
$N^\vee$ generated by $u_i$ and $N$ for $i=1,2,3.$  It can be
verified that all these three overlattices of $N$ satisfy the
conditions in Theorem \ref{urabe}. Therefore there are three reduced
sextic curves $D_1,D_2,D_3$ whose discriminantal groups are
$M_1^\vee/M_1,M_2^\vee/M_2$ and $M_3^\vee/M_3$ respectively. By the
algorithm in \cite{yang} we determine that all these three curves
have the same configuration. That is the union of an irreducible
quartic curve with an $A_3$ singularity and a conic such that the
conic meets the quartic at an $A_{15}$ point. Thus we have the
following theorem

\begin{thm}\label{triple}
The collection of three sextic curves $\{D_1,D_2,D_3\}$ forms a
Zariski triplet.
\end{thm}

These three curves are distinguished by the following two
conditions:

1) The line passing $A_{15}$ and $A_3$ is tangent to the quartic
component for $D_3,$ but not for $D_1$ and $D_2.$

2) For $D_1$ and $D_2,$ let $Q$ be the conic passing through the two
singularities together with their infinitely near points such that
the intersection number of the conic component and $Q$ at $A_{15}$
is at least $3.$ Then $Q \cap D_2 = \{ A_{15}, A_3 \}$ while $Q \cap
D_1$ contains four distinct points.

The equations of the curves are:

\medskip

$$D_1:(2x_0x_2+2x_1^2+x_2^2)[(3\lambda+2)x_1^4-4x_0^3x_2+\lambda x_0^2x_2^2
-4\lambda x_0^2x_1x_2-4x_0^2x_1^2$$
$$-2\lambda
x_0x_1x_2^2+(4\lambda+4)x_0x_1^2x_2-4\lambda x_0x_1^3+\lambda
x_1^2x_2^2]=0$$

for a generic parameter $\lambda$ and the $A_{15},A_3$ are located
at $(1:0:0)$ and $(0:0:1)$ respectively.

\medskip

$$D_2:(x_0x_2-x_1^2+\lambda x_2^2)
(x_0^3x_2-x_0^2x_1^2-5 \lambda x_1^4+10 \lambda x_0x_1^2x_2-4
\lambda x_0^2x_2^2)=0$$

for a generic parameter $\lambda$ and the $A_{15},A_3$ are located
at $(1:0:0)$ and $(0:0:1)$ respectively. the conic $x_0x_2-x_1^2$
intersects the sextic at $A_{15}$ and $A_3$ only.

\medskip

$$D_3:(x_0x_2+x_1^2-x_1x_2+x_2^2)(x_0^3x_2+x_0^2x_1^2+\lambda x_2^4-x_0^2x_1x_2+x_0^2x_2^2)=0,$$

and the $A_{15},A_3$ are located at $(1:0:0)$ and $(0:1:0)$
respectively. The line $x_2=0$ passes $A_{15},A_3$ and is tangent to
the conic component at $(1:0:0).$

\subsection{More than one Zariski pairs for the same combination of
singularities}

Our search shows several occasions where more than one Zariski pairs
pop up for the same combination of singularities. Here is a typical
example.

\begin{thm}\label{four_pairs}
There are four Zariski pairs of reduced sextic curves with
$3A_5+2A_1$ as their singularities. Their configurations are as
follows:

1) a quintic plus a line;

2) a quartic plus a conic;

3) two cubics;

4) a cubic plus a conic plus a line.

\end{thm}

Proof. Let $L$ be the root lattice of the Dynkin graph $3A_5+2A_1.$
Let $N={\mathbb Z}\lambda \oplus L.$ Denote the $17$ generators of
$L$ by $e_i (1 \leq i \leq 17)$ in a natural way such that
$e_ie_{i+1}=1$ for $1 \leq i \leq 4, 6 \leq i \leq 9$ and $11 \leq i
\leq 14.$ Let
$$u_1=(\lambda+\sum_{i=1}^5 ie_i+e_{16}+e_{17})/2,$$
$$u_2=\lambda/2+\sum_{i=1}^5 ie_i/6+\sum_{i=1}^5 i(e_{i+5}+e_{i+10})/3+(e_{16}+e_{17})/2,$$
$$v_1=\sum_{i=1}^5 i(e_i+e_{i+5})/2+(e_{16}+e_{17})/2,$$
$$v_2=\sum_{i=1}^5 i(e_i+e_{i+5})/6+\sum_{i=1}^5 ie_{i+10}/3 +(e_{16}+e_{17})/2,$$
$$w_1=(\lambda+\sum_{i=1}^5 i(e_i+e_{i+5}+e_{i+10}))/2,$$
$$w_2=\lambda/2+\sum_{i=1}^5 i(e_i+e_{i+5}+e_{i+10})/6,$$
$$s_1=(\sum_{i=1}^5 i(e_{i+5}+e_{i+10})+e_{16}+e_{17})/2,$$
$$t_1=(\lambda+\sum_{i=1}^5 ie_i+e_{16}+e_{17})/2,$$
$$s_2=\sum_{i=1}^5 i(e_i+e_{i+5})/6+\sum_{i=1}^5 ie_{i+10}/3+(e_{16}+e_{17})/2,$$
$$t_2=(\lambda+\sum_{i=1}^5 ie_i + e_{16}+e_{17})/2,$$

Let $A_i$ be the sublattice of $N^\vee$ generated by $u_i$ and $N$
for $i=1,2.$ Let $B_i$ be the sublattice of $N^\vee$ generated by
$v_i$ and $N$ for $i=1,2.$ Let $D_i$ be the sublattice of $N^\vee$
generated by $w_i$ and $N$ for $i=1,2.$ Let $G_i$ be the sublattice
of $N^\vee$ generated by $s_i,t_i$ and $N$ for $i=1,2.$

It can be verified that all these eight overlattices of $N$ satisfy
the conditions in Theorem \ref{urabe}. Therefore there are eight
reduced sextic curves $C_1,C_2,E_1,E_2,P_1,P_2,Q_1,Q_2$ whose
discriminantal groups are
$$A_1^\vee/A_1,A_2^\vee/A_2,B_1^\vee/B_1,B_2^\vee/B_2,
D_1^\vee/D_1,D_2^\vee/D_2,G_1^\vee/G_1,G_2^\vee/G_2$$ respectively.

By using the algorithm in \cite{yang} we determine the
configurations of these curves are:

1) a quintic plus a line for $C_1$ and $C_2;$

2) a quartic plus a conic for $E_1$ and $E_2;$

3) two cubics for $P_1$ and $P_2;$

4) a cubic plus a conic plus a line for $Q_1$ and $Q_2.$

It follows that the pairs $\{C_1,C_2\},\{E_1,E_2\},\{P_1,P_2\}$ and
$\{Q_1,Q_2\}$ are Zariski pairs. $\Box$

\section{Other Zariski pairs with different discriminantal groups}

Other Zariski pairs can be obtained using the method in the main
text. They are listed below. The notations in the tables are briefly
explained as follows.

The configuration of a sextic is described by the incidence table
with its entries to be the local components of a simple singularites
according to the convention in Figure 6.

For example,

\begin{tabular}{ccc}
 & 3 & 3  \\
 $A_{17}$ & I & II \\
 $A_1$ & I,II & \end{tabular}

means that the sextic has a nodal cubic and a smooth cubic as its
irreducible components and they intersects at an $A_{17}$ point.

\thicklines

\begin{center}
  \setlength{\unitlength}{0.08cm}
  \begin{picture}(120,45)

        \bezier{90}(10,20)(19,20)(28,7)
        \bezier{90}(10,20)(19,20)(28,33)
        \bezier{80}(90,20)(99,20)(108,7)
        \bezier{80}(90,20)(99,20)(108,33)
        \bezier{80}(90,20)(81,20)(72,7)
        \bezier{80}(90,20)(81,20)(72,33)

    \put(-5,20){\makebox(0,0){${\mathbf A}_{2n}:$}}
    \put(60,20){\makebox(0,0){${\mathbf A}_{2n-1}:$}}
    \put(20,4){\makebox(0,0){I}}
    \put(80,30){\makebox(0,0){I}}
    \put(80,4){\makebox(0,0){II}}

  \end{picture}
\end{center}

\begin{center}
  \setlength{\unitlength}{0.08cm}
  \begin{picture}(120,45)

        \put(17,20){\line(1,0){25}}
        \bezier{80}(30,20)(30,14)(40,8)
        \bezier{80}(30,20)(30,26)(40,32)
        \bezier{80}(30,20)(30,14)(20,8)
        \bezier{80}(30,20)(30,26)(20,32)

        \put(80,10){\line(1,0){30}}
        \bezier{80}(95,10)(95,22)(110,34)
        \bezier{80}(95,10)(95,22)(80,34)

    \put(-5,20){\makebox(0,0){${\mathbf D}_{2n}(n \geq 2):$}}
    \put(65,20){\makebox(0,0){${\mathbf D}_{2n+1}(n \geq 2):$}}
    \put(20,38){\makebox(0,0){I}}
    \put(40,38){\makebox(0,0){II}}
    \put(12,20){\makebox(0,0){III}}
    \put(80,4){\makebox(0,0){I}}
    \put(90,38){\makebox(0,0){II}}

  \end{picture}
\end{center}

\begin{center}
  \setlength{\unitlength}{0.08cm}
  \begin{picture}(120,40)

        \put(65,20){\line(1,0){25}}

    \bezier{80}(10,20)(20,20)(27,10)
        \bezier{80}(10,20)(20,20)(27,30)

        \bezier{80}(110,20)(120,20)(128,10)
        \bezier{80}(110,20)(120,20)(128,30)

        \bezier{80}(80,20)(73,20)(66,10)
        \bezier{80}(80,20)(73,20)(66,30)

    \put(-5,20){\makebox(0,0){${\mathbf E}_6:$}}
    \put(50,20){\makebox(0,0){${\mathbf E}_7:$}}
    \put(100,20){\makebox(0,0){${\mathbf E}_8:$}}
    \put(60,20){\makebox(0,0){I}}
    \put(80,4){\makebox(0,0){II}}
    \put(12,4){\makebox(0,0){I}}
    \put(116,4){\makebox(0,0){I}}

  \end{picture}
\end{center}

\begin{center}
{\bf Figure 6: Local components of simple singularities}
\end{center}

The overlattices are denoted by the generators over the root
lattice. The trivial lattice, i.e. the root lattice, is denote by
``-''. A generator of the overlattice is represented by its image in
$L^\vee/L$ where $L$ is the root lattice. As long as there is no
singularity of type $D_{2n},$ the discriminantal group of every
singularity is a cyclic group. Moreover, $({\mathbb
Z}\lambda)^\vee/{\mathbb Z}\lambda \cong {\mathbb Z}/2{\mathbb Z}$
is also cyclic. Hence a generator can be represented by a sequence
of numerals of which the first corresponds to its component in
$({\mathbb Z}\lambda)^\vee/{\mathbb Z}\lambda$.

Take the example of $D_7+A_{11}:$ the element $026$ stands for
$2u+6v$ where $u$ and $v$ are generators of $L(D_7)^\vee/L(D_7)$ and
$L(A_{11})^\vee/L(A_{11})$  respectively.

The last column in the table is a special curve distinguishing two
members of the Zariski pair. It is given by the degree of the curve
followed by a sequence of intersection numbers with the sextic curve
at singularities.

\pagebreak

\begin{center}
{\bf Table 1. Milnor number $18$}
\end{center}

\small

\setlength{\hoffset}{-20mm} \vsize =21.5truecm \hsize = 14.5 truecm
\begin{tabular}{|l|l|l|l|}
 \hline singularities & configuration & overlattices & special curve \\
\hline
 $3E_6$ & irreducible  & -,0111 & 2,$(E_6,4),(E_6,4),(E_6,4)$ \\
\hline
 $2E_6+A_5+A_1$ & irreducible  & -,01120 & 2,$(E_6,4),(E_6,4),(A_5,4)$ \\
\hline
 $E_6+A_{11}+A_1$ & irreducible  & -,0140 & 2,$(E_6,4),(A_{11},8)$ \\
\hline
 $E_6+A_8+A_2+2A_1$ & irreducible  & -,013100 & 2,$(E_6,4),(A_8,6),(A_2,2)$ \\
\hline
 $E_6+2A_5+2A_1$ &  \begin{tabular}{cccccc}
 & $E_6$ &  $A_5$ & $A_5$ & $A_1$ & $A_1$ \\
 5 & I & I,II & I & I & I \\
 1 & & & II & II & II \end{tabular}
 & 100311, 111211  &  2,$(E_6,4),(A_5,4),(A_5,4)$ \\
\hline
 $E_6+2A_5+2A_1$ &  \begin{tabular}{cccccc}
 & $E_6$ &  $A_5$ & $A_5$ & $A_1$ & $A_1$ \\
 4 & I & I & I & I & I \\
 2 & & II & II & II & II \end{tabular}
 & 003311, 011111  &
 2, $\begin{array}{l} (E_6,4),(A_5,2),(A_5,2), \\ (A_1,2),(A_1,2) \end{array}$ \\
\hline
 $D_7+A_{11}$ &  \begin{tabular}{ccc}
 & $D_7$ &  $A_{11}$ \\
 4 & II & II \\
 2 & I & I \end{tabular}
 & 026,013 &  2,$(D_7,6),(A_{11},6)$ \\
\hline
 $D_7+A_7+A_3+A_1$ &  \begin{tabular}{ccccc}
 & $D_7$ & $A_7$ & $A_3$ & $A_1$ \\
 4 & II & I & I &   \\
 1 & I &  &  II & I  \\
 1 &  & II &  & II \end{tabular}
 & \begin{tabular}{c} \{02420,10401\} \\ \{01210,10401\} \end{tabular}
  &  2, $(D_7,6),(A_7,4),(A_3,2)$\\
\hline
 $A_{17}+A_1$ &  irreducible
 & -,060 & 2,$(A_{17},12)$ \\
\hline
 $A_{17}+A_1$ &  \begin{tabular}{ccc}
 & $A_{17}$ &  $A_1$ \\
 3 & I & I,II \\
 3 & II & \end{tabular}
 & 190,130 &  1,$(A_{17},6)$ \\
\hline
 $A_{14}+A_2+2A_1$ &  irreducible
 & -,05100 & 2,$(A_{14},10),(A_2,2)$ \\
\hline
 $A_{11}+A_5+2A_1$ &  \begin{tabular}{ccccc}
 & $A_{11}$ & $A_5$ & $A_1$ & $A_1$ \\
 5 & I,II & I & I & I  \\
 1 &  & II & II & II \end{tabular}
 & 10311, 14111 &  2,$(A_{11},8),(A_5,4)$\\
\hline
 $A_{11}+A_5+2A_1$ &  \begin{tabular}{ccccc}
 & $A_{11}$ & $A_5$ & $A_1$ & $A_1$ \\
 4 & I & I,II & I & I  \\
 2 & II & & II & II \end{tabular}
 & 10311, 14111 &  2,$(A_{11},8),(A_5,4)$\\
\hline
 $A_{11}+A_5+2A_1$ &  \begin{tabular}{ccccc}
 & $A_{11}$ & $A_5$ & $A_1$ & $A_1$ \\
 3 & I & I &  I,II &   \\
 3 & II & II &  & I,II \end{tabular}
 & 10311, 14111 &  1,$(A_{11},4),(A_5,2)$\\
\hline
 $A_{11}+A_5+2A_1$ &  \begin{tabular}{ccccc}
 & $A_{11}$ & $A_5$ & $A_1$ & $A_1$ \\
 3 & I & I &  &   \\
 2 & II &  &  I & I  \\
 1 &  & II & II & II \end{tabular}
 & \begin{tabular}{c} \{06011,10311\} \\ \{02211,10311\} \end{tabular}
  &  2,$\begin{array}{l} (A_{11},4),(A_5,4), \\ (A_1,2),(A_1,2) \end{array}$\\
\hline
 $A_{11}+2A_2+3A_1$ &  \begin{tabular}{ccccccc}
 & $A_{11}$ & $A_2$ & $A_2$ & $A_1$ & $A_1$ & $A_1$ \\
 4 & I & I & I & I,II & I & I \\
 2 & II & &  &  & II & II \end{tabular}
 & 0600011,  0211011
  &  2, $\begin{array}{l} (A_{11},4),(A_2,2),(A_2,2), \\ (A_1,2),(A_1,2) \end{array}$ \\
\hline
 $A_{11}+2A_2+3A_1$ &  \begin{tabular}{ccccccc}
 & $A_{11}$ & $A_2$ & $A_2$ & $A_1$ & $A_1$ & $A_1$ \\
 3 & I & I & & I & I & I \\
 3 & II & & I & II & II & II \end{tabular}
 & 1600111, 1211111 &  2,$(A_{11},8),(A_2,2),(A_2,2)$\\
\hline
 $2A_9$ &  irreducible
 & -,024 & 2,$(A_9,4),(A_9,8)$ \\
\hline
\end{tabular}

\pagebreak

\begin{center}
{\bf Table 1 (cont.)}
\end{center}

\small

\setlength{\hoffset}{-20mm} \vsize =21.5truecm \hsize = 14.5 truecm
\begin{tabular}{|l|l|l|l|}
 \hline singularities & configuration & overlattices & special curve \\
\hline
 $2A_9$ &  \begin{tabular}{ccc}
 & $A_9$ &  $A_9$ \\
 5 & I,II & I \\
 1 &  & II \end{tabular}
 & 105,121 &  1,$(A_9,4),(A_9,2)$ \\
\hline
 $A_9+2A_4+A_1$ &  irreducible & -,02220 &
 2,$(A_9,4),(A_4,4),(A_4,4)$ \\
\hline
 $A_9+2A_4+A_1$ &  \begin{tabular}{ccccc}
 & $A_9$ &  $A_4$ & $A_4$ & $A_1$ \\
 5 & I & I & I & I,II \\
 1 & II & &  &  \end{tabular}
 & 15000, 11110 & 1,$(A_9,2),(A_4,2),(A_4,2)$ \\
\hline
 $2A_8+2A_1$ &  irreducible & -,03300 &
 2,$(A_8,6),(A_8,6)$ \\
\hline
 $A_8+A_5+A_2+3A_1$ &  \begin{tabular}{ccccccc}
 & $A_8$ &  $A_5$ & $A_2$ & $A_1$ & $A_1$ & $A_1$ \\
 5 & I & I & I & I,II & I & I\\
 1 &  & II &  & & II & II \end{tabular}
 & 1030011, 1311011 &  2,$(A_8,6),(A_5,4),(A_2,2)$ \\
\hline
 $2A_7+A_3+A_1$ &  \begin{tabular}{ccccc}
 & $A_7$ &  $A_7$ & $A_3$ & $A_1$ \\
 4 & I & I & I,II & I,II \\
 2 & II & II &  &  \end{tabular}
 & 01311, 12221 &
 2, $\begin{array}{l} (A_7,6),(A_7,2), \\ (A_3,2),(A_1,2) \end{array}$ \\
\hline
 $3A_5+3A_1$ &  \begin{tabular}{ccccccc}
 & $A_5$ &  $A_5$ & $A_5$ & $A_1$  & $A_1$ & $A_1$\\
 4 & I,II & I & I & I & I & \\
 1 & & II &  & II & & I \\
 1 & & & II & & II & II \end{tabular}
 & \begin{tabular}{c} \{0033011,1003101\} \\ \{0211011,1003101\} \end{tabular} &
 2, $\begin{array}{l} (A_5,4),(A_5,2),(A_5,2) \\ (A_1,2),(A_1,2) \end{array}$ \\
\hline
 $3A_5+3A_1$ &  \begin{tabular}{ccccccc}
 & $A_5$ &  $A_5$ & $A_5$ & $A_1$  & $A_1$ & $A_1$\\
 3 & I & I & I & I,II &  & \\
 2 & & II & II &  & I & I \\
 1 & II & &  & & II & II \end{tabular}
 & \begin{tabular}{c} \{0033011,1300011\} \\ \{0211011,1300011\} \end{tabular} &
 2, $\begin{array}{l} (A_5,4),(A_5,2),(A_5,2) \\ (A_1,2),(A_1,2) \end{array}$ \\
\hline
\end{tabular}

\begin{center}
{\bf Table 2. Milnor number $17$}
\end{center}

\small

\setlength{\hoffset}{-20mm} \vsize =21.5truecm \hsize = 14.5 truecm
\begin{tabular}{|l|l|l|l|}
 \hline singularities & configuration & overlattices & special curve \\
\hline
 $2E_6+A_5$ & irreducible  & -,0112 & 2,$(E_6,4),(E_6,4),(A_5,4)$ \\
\hline
 $2E_6+2A_2+A_1$ & irreducible  & -,01110 &
 2,$\begin{array}{l} (E_6,4),(E_6,4), \\(A_2,2),(A_2,2) \end{array}$ \\
\hline
 $E_6+A_{11}$ & irreducible  & -,014 & 2,$(E_6,4),(A_{11},8)$ \\
\hline
 $E_6+A_8+A_2+A_1$ & irreducible  & -,01310 & 2,$(E_6,4),(A_8,6),(A_2,2)$ \\
\hline
 $E_6+2A_5+A_1$ &  irreducible
 & -,01220 & 2,$(E_6,4),(A_5,4),(A_5,4)$ \\
\hline
 $E_6+A_5+2A_2+2A_1$ &  irreducible
 & -,0121100 & 2,$\begin{array}{l} (E_6,4),(A_5,4),\\(A_2,2),(A_2,2) \end{array}$ \\
\hline
 $E_6+A_5+2A_2+2A_1$ &  \begin{tabular}{ccccccc}
 & $E_6$ &  $A_5$ & $A_2$ & $A_2$ & $A_1$ & $A_1$ \\
 5 & I & I & I & I & I & I\\
 1 & & II & & & II & II \end{tabular}
 & 1030011,1111111  &
 2, $\begin{array}{l} (E_6,4),(A_5,4), \\ (A_2,2),(A_2,2) \end{array}$ \\
\hline
 $D_7+A_7+A_3$ &  \begin{tabular}{cccc}
 & $D_7$ & $A_7$ & $A_3$ \\
 4 & II & I & I \\
 2 & I & II & II \end{tabular}
 & 0242,0121 & 2,$(D_7,6),(A_7,4),(A_3,2)$ \\
\hline
 $D_7+3A_3+A_1$ &  \begin{tabular}{cccccc}
 & $D_7$ &  $A_3$ & $A_3$ & $A_3$ & $A_1$\\
 4 & II & I  & I & I & \\
 1 & & & II & II & I \\
 1 & I & II & & & II \end{tabular}
 & \begin{tabular}{c} \{022220,100221\} \\ \{011110,100221\} \end{tabular} &
 2, $\begin{array}{l} (D_7,6),(A_3,2) \\ (A_3,2),(A_3,2) \end{array}$ \\
\hline
 $2D_5+A_7$ &  \begin{tabular}{cccc}
 & $D_5$ & $D_5$ & $A_7$ \\
 4 & II & II & I \\
 2 & I & I & II \end{tabular}
 & 0224,0112 & 2,$(D_5,4),(D_5,4),(A_7,4)$ \\
\hline
\end{tabular}

\pagebreak

\begin{center}
{\bf Table 2 (cont.)}
\end{center}

\small

\setlength{\hoffset}{-20mm} \vsize =21.5truecm \hsize = 14.5 truecm
\begin{tabular}{|l|l|l|l|}
 \hline singularities & configuration & overlattices & special curve \\
\hline
 $2D_5+2A_3+A_1$ &  \begin{tabular}{cccccc}
 & $D_5$ & $D_5$ & $A_3$ & $A_3$ & $A_1$ \\
 4 & II & II & I & I &\\
 1 & I & & II & & I \\
 1 & & I & &  II & II \end{tabular}
 & \begin{tabular}{c} \{022220,102021\} \\ \{011110,102021\} \end{tabular} &
 2,$\begin{array}{l} (D_5,4),(D_5,4) \\ (A_3,2),(A_3,2) \end{array}$ \\
\hline
 $D_5+A_{11}+A_1$ &  \begin{tabular}{cccc}
 & $D_5$ & $A_{11}$ & $A_1$ \\
 4 & II & I & I,II \\
 2 & I & II &  \end{tabular}
 & 0260,0131 & 2,$(D_5,4),(A_{11},6),(A_1,2)$ \\
\hline
 $D_5+A_7+A_3+2A_1$ &  \begin{tabular}{cccccc}
 & $D_5$ & $A_7$ & $A_3$ & $A_1$ & $A_1$\\
 4 & II & I & I & I,II &\\
 1 & I & & II & & I \\
 1 & & II & & & II  \end{tabular}
 & \begin{tabular}{c} \{024200,104001\} \\ \{012110,104001\} \end{tabular} &
 2,$\begin{array}{l} (D_5,4),(A_7,4) \\ (A_3,2),(A_1,2) \end{array}$ \\
\hline
 $A_{17}$ &  irreducible
 & -,06 & 2,$(A_{17},12)$ \\
\hline
 $A_{17}$ &  \begin{tabular}{cc}
 & $A_{17}$ \\
 3 & I \\
 3 & II \end{tabular}
 & 19,13 &  1,$(A_{17},6)$ \\
\hline
 $A_{15}+2A_1$ &  \begin{tabular}{cccc}
 & $A_{15}$ & $A_1$ & $A_1$ \\
 4 & I & I,II & I,II \\
 2 & II & & \end{tabular}
 & 0800,0411 & 2, $(A_{15},8),(A_2,2),(A_1,2)$ \\
\hline
 $A_{14}+A_2+A_1$ &  irreducible
 & -,0510 & 2,$(A_{14},10),(A_2,2)$ \\
\hline
 $A_{11}+A_5+A_1$ &  irreducible
 & -,0420 & 2,$(A_{11},8),(A_5,4)$ \\
\hline
 $A_{11}+A_5+A_1$ &  \begin{tabular}{cccc}
 & $A_{11}$ & $A_5$ & $A_1$ \\
 3 & I & I & I,II\\
 3 & II & II & \end{tabular}
 & 1630,1210 & 1,$(A_{11},4),(A_5,2)$ \\
\hline
 $A_{11}+2A_3$ &  \begin{tabular}{cccc}
 & $A_{11}$ & $A_3$ & $A_3$ \\
 4 & I & I & I,II \\
 2 & II & II &  \end{tabular}
 & 0620,0312 & 2,$(A_{11},6),(A_3,4),(A_3,2)$ \\
\hline
 $A_{11}+2A_2+2A_1$ & irreducible
 & -,041100 &  2,$(A_{11},8),(A_2,2),(A_2,2)$ \\
\hline
 $A_{11}+2A_2+2A_1$ &  \begin{tabular}{cccccc}
 & $A_{11}$ & $A_2$ & $A_2$ & $A_1$ & $A_1$  \\
 4 & I & I & I &  I & I \\
 2 & II & &  & II & II \end{tabular}
 & 060011,021111 & 2,$\begin{array}{l} (A_{11},4),(A_2,2),(A_2,2) \\ (A_1,2),(A_1,2) \end{array}$ \\
\hline
 $A_9+2A_4$ &  irreducible & -,0222 &
 2,$(A_9,4),(A_4,4),(A_4,4)$ \\
\hline
 $A_9+2A_4$ &  \begin{tabular}{cccc}
 & $A_9$ &  $A_4$ & $A_4$  \\
 5 & I & I & I \\
 1 & I & &  \end{tabular}
 & 1500,1111 & 1,$(A_9,2),(A_4,2),(A_4,2)$ \\
\hline
 $2A_8+A_1$ &  irreducible & -,0330 &
 2,$(A_8,6),(A_8,6)$ \\
\hline
 $A_8+A_5+A_2+2A_1$ &  irreducible & -,032100 &
 2,$(A_8,6),(A_5,4),(A_2,2)$ \\
\hline
 $A_8+A_5+A_2+2A_1$ &  \begin{tabular}{cccccc}
 & $A_8$ &  $A_5$ & $A_2$ & $A_1$ & $A_1$  \\
 5 & I & I & I & I & I\\
 1 & & II & & II & II \end{tabular}
 & 103011,131111 & 2,$(A_8,6),(A_5,4),(A_2,2)$ \\
\hline
 $A_8+3A_2+3A_1$ & irreducible
 & -,03111000 & 2,$\begin{array}{l} (A_8,6),(A_2,2) \\ (A_2,2),(A_2,2) \end{array}$
\\
\hline
 $2A_7+3A_1$ &  \begin{tabular}{cccccc}
 & $A_7$ &  $A_7$ & $A_1$ & $A_1$ & $A_1$ \\
 4 & I & I & I,II & I,II &  \\
 1 & II &  &  &  & I \\
 1 &  & II & & & II \end{tabular}
 & \begin{tabular}{c} \{044000,104001\} \\ \{022110,104001\} \end{tabular} &
 2,$\begin{array}{l} (A_7,4),A_7,4) \\ (A_1,2),(A_1,2) \end{array}$ \\
\hline
\end{tabular}

\pagebreak

\begin{center}
{\bf Table 2 (cont.)}
\end{center}

\small

\setlength{\hoffset}{-20mm} \vsize =21.5truecm \hsize = 14.5 truecm
\begin{tabular}{|l|l|l|l|}
 \hline singularities & configuration & overlattices & special curve \\
\hline
 $A_7+3A_3+A_1$ &  \begin{tabular}{cccccc}
 & $A_7$ &  $A_3$ & $A_3$ & $A_3$ & $A_1$ \\
 4 & I & I,II & I & I &  \\
 1 & &  & II & II & I \\
 1 & II &  & & & II \end{tabular}
 & \begin{tabular}{c} \{040220,100221\} \\ \{021120,102201\} \end{tabular} &
 2,$\begin{array}{l} (A_7,4),A_3,4) \\ (A_3,2),(A_3,2) \end{array}$ \\
\hline
 $2A_5+2A_2+3A_1$ &  \begin{tabular}{cccccccc}
 & $A_5$ &  $A_5$ & $A_2$ & $A_2$ & $A_1$ & $A_1$ & $A_1$  \\
 5 & I,II & I & I & I & I,II & I & I \\
 1 & & II & & & & II & II \end{tabular}
 & 10300011,11211011 & 2,$\begin{array}{l} (A_5,4),A_5,4) \\ (A_2,2),(A_2,2) \end{array}$ \\
\hline
 $2A_5+2A_2+3A_1$ &  \begin{tabular}{cccccccc}
 & $A_5$ &  $A_5$ & $A_2$ & $A_2$ & $A_1$ & $A_1$ & $A_1$  \\
 4 & I & I & I & I & I,II & I & I \\
 2 & II & II & & & & II & II \end{tabular}
 & 03300011,01111011 & 2,$\begin{array}{l} (A_5,2),(A_5,2),(A_2,2) \\ (A_2,2),(A_1,2),(A_1,2) \end{array}$ \\
\hline
 $2A_5+2A_2+3A_1$ &  \begin{tabular}{cccccccc}
 & $A_5$ &  $A_5$ & $A_2$ & $A_2$ & $A_1$ & $A_1$ & $A_1$  \\
 3 & I & I & I &  & I & I & I \\
 3 & II & II & & I & II & II & II \end{tabular}
 & 13300111,11111111 & 2,$\begin{array}{l} (A_5,4),A_5,4) \\ (A_2,2),(A_2,2) \end{array}$ \\
\hline
 $2A_5+2A_2+3A_1$ &  \begin{tabular}{cccccccc}
 & $A_5$ &  $A_5$ & $A_2$ & $A_2$ & $A_1$ & $A_1$ & $A_1$  \\
 4 & I & I & I & I &  & I & I \\
 1 & II & & & & I & II &  \\
 1 & & II & & & II & & II \end{tabular}
 & \begin{tabular}{c} \{03300011,10300101\} \\ \{01111011,10300101\} \end{tabular}
 & 2,$\begin{array}{l} (A_5,2),A_5,2),(A_2,2) \\ (A_2,2),(A_1,2),(A_1,2) \end{array}$ \\
\hline
 $4A_4+A_1$ & irreducible
 & -,011220 & 2,$\begin{array}{l} (A_4,4),A_4,4) \\ (A_4,2),(A_4,2) \end{array}$ \\
\hline
\end{tabular}

\begin{center}
{\bf Table 3. Milnor number $16$}
\end{center}

\small

\setlength{\hoffset}{-20mm} \vsize =21.5truecm \hsize = 14.5 truecm
\begin{tabular}{|l|l|l|l|}
 \hline singularities & configuration & overlattices & special curve \\
\hline
 $2E_6+2A_2$ & irreducible  & -,0111 &
 2,$\begin{array}{l} (E_6,4),(E_6,4),\\ (A_2,2),(A_2,2) \end{array}$ \\
\hline
 $E_6+A_8+A_2$ & irreducible  & -,0131 & 2,$(E_6,4),(A_8,6),(A_2,2)$ \\
\hline
 $E_6+2A_5$ &  irreducible
 & -,0122 & 2,$(E_6,4),(A_5,4),(A_5,4)$ \\
\hline
 $E_6+A_5+2A_2+A_1$ &  irreducible
 & -,012110 & 2,$\begin{array}{l} (E_6,4),(A_5,4),\\ (A_2,2),(A_2,2) \end{array}$ \\
\hline
 $E_6+4A_2+2A_1$ &  irreducible
 & -,01111100 & 2,$\begin{array}{l} (E_6,4),(A_2,2),(A_2,2),\\ (A_2,2),(A_2,2) \end{array}$ \\
\hline
 $D_7+3A_3$ &  \begin{tabular}{ccccc}
 & $D_7$ &  $A_3$ & $A_3$ & $A_3$ \\
 4 & II & I & I & I \\
 2 & I & II & II & II \\
 \end{tabular}
 & 02222,01111 & 2,$\begin{array}{l} (D_7,6),(A_3,2), \\ (A_3,2),(A_3,2) \end{array}$ \\
\hline
 $2D_5+2A_3$ &  \begin{tabular}{ccccc}
 & $D_5$ &  $D_5$ & $A_3$ & $A_3$ \\
 4 & II & II & I & I \\
 2 & I & I & II & II \\
 \end{tabular}
 & 02222,01111 & 2,$\begin{array}{l} (D_5,4),(D_5,4), \\ (A_3,2),(A_3,2) \end{array}$ \\
\hline
 $D_5+A_7+A_3+A_1$ &  \begin{tabular}{ccccc}
 & $D_5$ &  $A_7$ & $A_3$ & $A_1$ \\
 4 & II & I & I & I,II\\
 2 & I & II & II & \\
 \end{tabular}
 & 02420,01211 &  2, $\begin{array}{l} (D_5,4),(A_7,4), \\ (A_3,2),(A_1,2) \end{array}$ \\
\hline
 $D_5+3A_3+2A_1$ &  \begin{tabular}{ccccccc}
 & $D_5$ & $A_3$ & $A_3$ & $A_3$ & $A_1$ & $A_1$\\
 4 & II & I & I & I & & I,II \\
 1 & I & II & & & I & \\
 1 & &  & II & II & II &  \end{tabular}
 & \begin{tabular}{c} \{0222200,1002201\} \\ \{0111110,1002201\} \end{tabular} &
 2,$\begin{array}{l} (D_5,4),(A_3,2),(A_3,2) \\ (A_3,2),(A_1,2) \end{array}$ \\
\hline
\end{tabular}

\pagebreak

\begin{center}
{\bf Table 3 (cont.)}
\end{center}

\small

\setlength{\hoffset}{-20mm} \vsize =21.5truecm \hsize = 14.5 truecm
\begin{tabular}{|l|l|l|l|}
 \hline singularities & configuration & overlattices & special curve \\
\hline
 $A_{14}+A_2$ &  irreducible
 & -,051 & 2,$(A_{14},10),(A_2,2)$ \\
\hline
 $A_{11}+A_5$ &  irreducible
 & -,042 & 2,$(A_{11},8),(A_5,4)$ \\
\hline
 $A_{11}+A_5$ &  \begin{tabular}{ccc}
 & $A_{11}$ & $A_5$ \\
 3 & I & I \\
 3 & II & II \end{tabular}
 & 163,121 & 1,$(A_{11},4),(A_5,2)$ \\
\hline
 $A_{11}+A_3+2A_1$ &  \begin{tabular}{ccccc}
 & $A_{11}$ & $A_3$ & $A_1$ & $A_1$ \\
 4 & I & I  & I,II & I,II \\
 2 & II & II & & \end{tabular}
 & 06200,03111 &  2, $\begin{array}{l} (A_{11},6),(A_3,2), \\ (A_1,2),(A_1,2) \end{array}$ \\
\hline
 $A_{11}+2A_2+A_1$ & irreducible
 & -,04110 &  2,$(A_{11},8),(A_2,2),(A_2,2)$ \\
\hline
 $2A_8$ &  irreducible & -,033 &
 2,$(A_8,6),(A_8,6)$ \\
\hline
 $A_8+A_5+A_2+A_1$ & irreducible
 & -,03210 & 2,$(A_8,6),(A_5,4),(A_2,2)$ \\
\hline
 $A_8+3A_2+2A_1$ & irreducible
 & -,0311100 & 2,$\begin{array}{l} (A_8,6),(A_2,2), \\ (A_2,2),(A_2,2) \end{array}$ \\
\hline
 $2A_7+2A_1$ &  \begin{tabular}{ccccc}
 & $A_7$ &  $A_7$ & $A_1$ & $A_1$\\
 4 & I & I & I,II & I,II \\
 2 & II & II &  &  \end{tabular}
 & 04400,02211 & 2, $\begin{array}{l} (A_7,4),(A_7,4), \\ (A_1,2),(A_1,2) \end{array}$ \\
\hline
 $A_7+3A_3$ &  \begin{tabular}{ccccc}
 & $A_7$ &  $A_3$ & $A_3$ & $A_3$\\
 4 & I & I,II & I & I \\
 2 & II & & II & II \end{tabular}
 & 04022,02112 & 2, $\begin{array}{l} (A_7,4),(A_3,4), \\ (A_3,2),(A_3,2) \end{array}$ \\
\hline
 $A_7+2A_3+3A_1$ &  \begin{tabular}{ccccccc}
 & $A_7$ & $A_3$ & $A_3$ & $A_1$ & $A_1$ & $A_1$\\
 4 & I & I & I & I,II & I,II & \\
 1 & II & & & & & I\\
 1 & & II & II &  &  & II \end{tabular}
 & \begin{tabular}{c} \{0422000,1022001\} \\ \{0211110,1022001\} \end{tabular} &
 2,$\begin{array}{l} (A_7,4),(A_3,2),(A_3,2) \\ (A_1,2),(A_1,2) \end{array}$ \\
\hline
 $3A_5+A_1$ & irreducible
 & -,02220 & 2,$(A_5,4),(A_5,4),(A_5,4)$ \\
\hline
 $3A_5+A_1$ &  \begin{tabular}{ccccc}
 & $A_5$ &  $A_5$ & $A_5$ & $A_1$\\
 3 & I & I & I & I,II \\
 3 & II & II & II & \end{tabular}
 & 13330,11110 &  1,$(A_5,2),(A_5,2),(A_5,2)$\\
\hline
 $2A_5+2A_2+2A_1$ & irreducible
 & -,0221100 & 2,$\begin{array}{l} (A_5,4),(A_5,4), \\ (A_2,2),(A_2,2) \end{array}$ \\
\hline
 $2A_5+2A_2+2A_1$ &  \begin{tabular}{ccccccc}
 & $A_5$ &  $A_5$ & $A_2$ & $A_2$ & $A_1$ & $A_1$\\
 5 & I,II & I & I & I & I & I\\
 1 & & II & & & II & II \end{tabular}
 & 1030011,1121111 &  2, $(A_5,6),(A_5,6)$ \\
\hline
 $2A_5+2A_2+2A_1$ &  \begin{tabular}{ccccccc}
 & $A_5$ &  $A_5$ & $A_2$ & $A_2$ & $A_1$ & $A_1$\\
 4 & I & I & I & I & I & I\\
 2 & II & II & & & II & II \end{tabular}
 & 0330011,0111111 &  2, $\begin{array}{l} (A_5,2),(A_5,2),(A_2,2) \\ (A_2,2),(A_1,2),(A_1,2) \end{array}$ \\
\hline
 $A_5+4A_2+3A_1$ & irreducible
 & -,021111000 & 2,$\begin{array}{l} (A_5,4),(A_2,2),(A_2,2), \\ (A_2,2),(A_2,2) \end{array}$ \\
\hline
 $A_5+4A_2+3A_1$ &  \begin{tabular}{ccccccccc}
 & $A_5$ &  $A_2$ & $A_2$ & $A_2$ & $A_2$ & $A_1$ & $A_1$ & $A_1$\\
 5 & I & I & I & I & I & I,II & I & I\\
 1 & II & & & & & & II & II \end{tabular}
 & 130000011,111111011 &
 2, $\begin{array}{l} (A_5,4),(A_2,2),(A_2,2), \\ (A_2,2),(A_2,2) \end{array}$  \\
\hline
 $4A_4$ & irreducible
 & -,02211 & 2,$\begin{array}{l} (A_4,4),(A_4,4), \\ (A_4,2),(A_4,2) \end{array}$ \\
\hline
 $5A_3+A_1$ &  \begin{tabular}{ccccccc}
 & $A_3$ & $A_3$ & $A_3$ & $A_3$ & $A_3$ & $A_1$\\
 4 & I,II & I & I & I,II &  & \\
 1 &  & II & II & & & I\\
 1 & &  &  & II & II & II \end{tabular}
 & \begin{tabular}{c} \{0022220,1000221\} \\ \{0111120,1002201\} \end{tabular} &
 2,$\begin{array}{l} (A_3,4),(A_3,2),(A_3,2) \\ (A_3,2),(A_3,2) \end{array}$ \\
\hline
\end{tabular}

\pagebreak

\begin{center}
{\bf Table 4. Milnor number $15$ or less}
\end{center}

\small

\setlength{\hoffset}{-20mm} \vsize =21.5truecm \hsize = 14.5 truecm
\begin{tabular}{|l|l|l|l|}
 \hline singularities & configuration & overlattices & special curve \\
\hline
 $E_6+A_5+2A_2$ &  irreducible
 & -,01211 & 2,$\begin{array}{l} (E_6,4),(A_5,4),\\ (A_2,2),(A_2,2) \end{array}$ \\
\hline
 $E_6+4A_2+A_1$ &  irreducible
 & -,0111110 & 2,$\begin{array}{l} (E_6,4),(A_2,2),(A_2,2),\\ (A_2,2),(A_2,2) \end{array}$ \\
\hline
 $D_5+3A_3+A_1$ &  \begin{tabular}{cccccc}
 & $D_5$ &  $A_3$ & $A_3$ & $A_3$ & $A_1$ \\
 4 & II & I & I & I & I,II\\
 2 & I & II & II & II & \\
 \end{tabular}
 & 022220,011111 &  2, $\begin{array}{l} (D_5,4),(A_3,2),(A_3,2) \\ (A_3,2),(A_1,2) \end{array}$ \\
\hline
 $A_{11}+2A_2$ & irreducible
 & -,0411 &  2,$(A_{11},8),(A_2,2),(A_2,2)$ \\
\hline
 $A_8+A_5+A_2$ & irreducible
 & -,0321 & 2,$(A_8,6),(A_5,4),(A_2,2)$ \\
\hline
 $A_8+3A_2+A_1$ & irreducible
 & -,031110 & 2,$\begin{array}{l} (A_8,6),(A_2,2), \\ (A_2,2),(A_2,2) \end{array}$ \\
\hline
 $A_7+2A_3+2A_1$ &  \begin{tabular}{cccccc}
 & $A_7$ &  $A_3$ & $A_3$ & $A_1$ & $A_1$ \\
 4 & I & I & I & I,II & I,II \\
 2 & II & II & II & & \end{tabular}
 & 042200,021111 & 2, $\begin{array}{l} (A_7,4),(A_3,2),(A_3,2) \\ (A_1,2),(A_1,2) \end{array}$ \\
\hline
 $3A_5$ & irreducible
 & -,0222 & 2,$(A_5,4),(A_5,4),(A_5,4)$ \\
\hline
 $3A_5$ &  \begin{tabular}{cccc}
 & $A_5$ &  $A_5$ & $A_5$ \\
 3 & I & I & I \\
 3 & II & II & II\end{tabular}
 & 1333,1111 &  1,$(A_5,2),(A_5,2),(A_5,2)$\\
\hline
 $2A_5+2A_2+A_1$ & irreducible
 & -,022110 & 2,$\begin{array}{l} (A_5,4),(A_5,4) \\ (A_2,2),(A_2,2) \end{array}$ \\
\hline
 $A_5+4A_2+2A_1$ & irreducible
 & -,0211110 & 2,$\begin{array}{l} (A_5,4),(A_2,2),(A_2,2) \\ (A_2,2),(A_2,2) \end{array}$ \\
\hline
 $A_5+4A_2+2A_1$ &  \begin{tabular}{cccccccc}
 & $A_5$ &  $A_2$ & $A_2$ & $A_2$ & $A_2$ & $A_1$ & $A_1$ \\
 5 & I & I & I & I & I & I & I\\
 1 & II & & & & & II & II \end{tabular}
 & 13000011,11111111 &
 2, $\begin{array}{l} (A_5,4),(A_2,2),(A_2,2), \\ (A_2,2),(A_2,2) \end{array}$  \\
\hline
\end{tabular}

\pagebreak

\begin{center}
{\bf Table 4 (cont.)}
\end{center}

\small

\setlength{\hoffset}{-20mm} \vsize =21.5truecm \hsize = 14.5 truecm
\begin{tabular}{|l|l|l|l|}
 \hline singularities & configuration & overlattices & special curve \\
\hline
 $5A_3$ &  \begin{tabular}{cccccc}
 & $A_3$ &  $A_3$ & $A_3$ & $A_3$ & $A_3$ \\
 4 & I,II & I  & I & I & I  \\
 2 & & II &  II & II & II \end{tabular}
 & 002222,011112 &
 2,$\begin{array}{l} (A_3,4),(A_3,2),(A_3,2), \\ (A_3,2),(A_3,2) \end{array}$ \\
\hline
 $4A_3+3A_1$ &  \begin{tabular}{cccccccc}
 & $A_3$ &  $A_3$ & $A_3$ & $A_3$ & $A_1$ & $A_1$ & $A_1$ \\
 4 & I & I  & I & I & I,II & I,II &  \\
 1 & II & II  &  &  &  & & I\\
 1 & &  &  II & II & & & II \end{tabular}
 & $\begin{array}{l} \{02222000,10022001\} \\ \{01111011,10022100\} \end{array}$ &
 2,$\begin{array}{l} (A_3,2),(A_3,2),(A_3,2), \\ (A_3,2),(A_1,2),(A_1,2) \end{array}$ \\
\hline
 $6A_2+3A_1$ &  irreducible & -,0111111000 &
 2,$\begin{array}{l} (A_2,2),(A_2,2),(A_2,2), \\
 (A_2,2),(A_2,2),(A_2,2) \end{array}$ \\
\hline
 $E_6+4A_2$ &  irreducible
 & -,011111 & 2,$\begin{array}{l} (E_6,4),(A_2,2),(A_2,2),\\ (A_2,2),(A_2,2) \end{array}$ \\
\hline
 $A_8+3A_2$ & irreducible
 & -,03111 & 2,$\begin{array}{l} (A_8,6),(A_2,2), \\ (A_2,2),(A_2,2) \end{array}$ \\
\hline
 $2A_5+2A_2$ & irreducible
 & -,02211 & 2,$\begin{array}{l} (A_5,4),(A_5,4), \\ (A_2,2),(A_2,2) \end{array}$ \\
\hline
 $A_5+4A_2+A_1$ & irreducible
 & -,0211110 & 2,$\begin{array}{l} (A_5,4),(A_2,2),(A_2,2),\\ (A_2,2),(A_2,2) \end{array}$ \\
\hline
 $4A_3+2A_1$ &  \begin{tabular}{ccccccc}
 & $A_3$ &  $A_3$ & $A_3$ & $A_3$ & $A_1$ & $A_1$  \\
 4 & I & I  & I & I & I,II & I,II  \\
 2 & II & II &  II & II & & \end{tabular}
 & $0222200,0111111$ &
 2,$\begin{array}{l} (A_3,2),(A_3,2),(A_3,2), \\ (A_3,2),(A_1,2),(A_1,2) \end{array}$ \\
\hline
 $6A_2+2A_1$ &  irreducible & -,011111100 &
 2,$\begin{array}{l} (A_2,2),(A_2,2),(A_2,2), \\
 (A_2,2),(A_2,2),(A_2,2) \end{array}$ \\
\hline
 $A_5+4A_2$ & irreducible
 & -,021111 & 2,$\begin{array}{l} (A_5,4),(A_2,2),(A_2,2),\\ (A_2,2),(A_2,2) \end{array}$ \\
\hline
 $6A_2+A_1$ &  irreducible & -,01111110 &
 2,$\begin{array}{l} (A_2,2),(A_2,2),(A_2,2), \\
 (A_2,2),(A_2,2),(A_2,2) \end{array}$ \\
\hline
\end{tabular}

\pagebreak

\large


\begin{thebibliography}{99}

\bibitem{Bartolo1} E. Artal-Bartolo, {\em Sur les couples de
Zariski,} J. Algebraic Geom. {\bf 3} (1994), 223-247.

\bibitem{Bartolo2} E. Artal-Bartolo, H. Tokunaga, {\em Zariski pairs of index 19 and Mordell-Weil groups of K3
surfaces,} Proc. London Math. Soc. (3), {\bf 80} (2000), no.1,
127-144.

\bibitem{BPV}
W. P. Barth, K. Hulek, C. A. M. Peters, and A. Van de Ven, {\em
Compact complex surfaces (2nd ed.),} Springer-Verlag, Berlin, 2004.

\bibitem{kplets}
A. Degtyarev, {\em Zariski k-plets via dessins d'enfants,} preprint,
arxiv.math/0710.0279v2,2008.

\bibitem{Namba}
M.Namba and H.Tsuchihashi, {\em On the Fundamental Groups of Galois
Covering Spaces of the Projective Plane,} Geometriae Dedicata {\bf
105} (2004) , 85-105.

\bibitem{nik}
V.V.Nikulin, {\em Integral symmetric bilinear forms and some of
their applications,} Math. USSR Izvestiya {\bf 14} (1980), 103-167.

\bibitem{Shimada09}
I. Shimada, {\em Lattice Zariski k-ples of plane sextic curves and
Z-splitting curves for double plance sextics,}  preprint,
arXiv:math/0903.3308v1,2009.


\bibitem{Shimada1}
I. Shimada, {\em Non-homeomorphic conjugate complex varieties,}
preprint, arXiv:math/0701115,2009.


\bibitem{Shimada2}
I. Shimada, {\em On Arithmetic Zariski Pairs in degree $6,$}
preprint,arxiv/math/0611596, to appear in Adv. Geom.

\bibitem{sextic}
T.Urabe, {\em Combinations of rational singularities on plane sextic
curves with the sum of Milnor numbers less than sixteen,} Banach
Center Publ. {\bf 20} (1988), 429-456.

\bibitem{ams}
T.Urabe, {\em Dynkin graphs and combinations of singularities on
plane sextic curves,} in Singularities, Proc., Univ. Iowa 1986 (R.
Randell, ed.), Contemporary Math. {\bf 90}, Amer. Math.
Soc.,Province,Rhode Island, 1989, 295-316.

\bibitem{yang}
J. Yang, {\em Sextic curves with simple singularities,} Tohoku Math.
J., {\bf 48} (2) (1996), 203-227.

\bibitem{zariski1}
O. Zariski, {\em On the Problem of Existence of Algebraic Functions
of Two Variables Possessing a Given Branch Curve,} Amer. J. Math.,
{\bf 51} (2) (1929), 305-328.
\end{thebibliography}
\end{document}